\documentclass[11pt]{amsart}
\usepackage{amsmath}
\usepackage{amsxtra}
\usepackage{amscd}
\usepackage{amsthm}
\usepackage{amsfonts}
\usepackage{amssymb}
\usepackage{eucal}
\usepackage[all]{xy}
\begin{document}

\newcommand{\g}{{\mathfrak{g}}}
\newcommand{\slt}{\mathfrak{sl}_2}
\newcommand{\slth}{\widehat{\mathfrak{sl}}_2}
\newcommand{\slnh}{\widehat{\mathfrak{sl}}_n}
\newcommand{\sltr}{\mathfrak{sl}_3}
\newcommand{\sln}{\mathfrak{sl}_n}
\newcommand{\slnm}{\mathfrak{sl}_{n-1}}
\newcommand{\nb}{\mathbf{n}}
\newcommand{\nub}{{\boldsymbol \nu}}
\newcommand{\mb}{\mathbf{m}}
\newcommand{\Nb}{\mathbf{N}}
\newcommand{\kb}{\mathbf{k}}
\newcommand{\mub}{{\boldsymbol \lambda}}
\newcommand{\sbld}{\mathbf{m}}
\newcommand{\sbf}{\mathbf{s}}

\newcommand{\nn}{\nonumber}
\newcommand{\bea}{\begin{eqnarray}}
\newcommand{\ena}{\end{eqnarray}}
\newcommand{\be}{\begin{eqnarray*}}
\newcommand{\en}{\end{eqnarray*}}

\newcommand{\res}{{\mathop{\rm res}}}
\newcommand{\id}{{\rm id}}
\newcommand{\ch}{{\rm ch}}
\newcommand{\End}{\mathop{{\rm End}}}
\newcommand{\tr}{\mathop{{\rm tr}}}
\newcommand{\gr}{\mathop{{\rm gr}}}
\newcommand{\bra}[1]{\langle #1 |}        
\newcommand{\ket}[1]{{| #1 \rangle}}      
\newcommand{\br}[1]{{\langle #1 \rangle}}  
\newcommand{\qbin}[2]
{{ 
\left[
\begin{matrix}{\displaystyle #1}\\
{\displaystyle #2}\end{matrix}
\right]
}}
\newcommand{\ord}{{\mathop{\rm ord}}}
\newcommand{\xt}{\tilde{x}}
\newcommand{\bep}{\bar{\epsilon}}

\numberwithin{equation}{section}
\newtheorem{thm}{Theorem}[section]
\newtheorem{prop}[thm]{Proposition}
\newtheorem{lem}[thm]{Lemma}
\newtheorem{cor}[thm]{Corollary}
\newtheorem{rem}{Remark}
\newtheorem{dfn}{Definition}

\newcommand{\gc}{{\g[t]}}
\newcommand{\zz}{{\mathcal Z}}
\newcommand{\fusn}{\circledast}
\newcommand{\wg}{\widehat{\g}}
\newcommand{\wt}{\mathop{\rm wt}}
\newcommand{\ver}{{\mathcal V}^{(k)}}
\newcommand{\usl}{U(\slt[t])}
\newcommand{\ide}{{\mathcal A}}
\newcommand{\simk}{{P_+^{(k)}}}
\newcommand{\ant}{S}
\newcommand{\hkr}{{\theta^\vee}}
\newcommand{\fil}{{\mathcal F}_\zz}
\newcommand{\gn}{{\mathfrak n}}
\newcommand{\hh}{{\mathfrak h}}
\newcommand{\ann}{{\rm Ann}}
\newcommand{\gX}{{\mathfrak X}}
\newcommand{\gY}{{\mathfrak Y}}
\newcommand{\ik}{{\mathfrak I}^{(k)}}
\newcommand{\tx}{\tilde{X}}
\newcommand{\txx}{\tilde{\gX}}
\newcommand{\fk}{{\mathcal K}}
\newcommand{\ga}{{\mathfrak a}}
\newcommand{\cpu}{U_0(\gc)}
\newcommand{\Ikla}{I^{(k)}_{\lambda}}
\newcommand{\Ila}{I_{\lambda}}
\newcommand{\Lkla}{L^{(k)}_{\lambda}}
\newcommand{\vkla}{v^{(k)}_{\lambda}}
\newcommand{\Lkinf}{L^{(k)}_{\lambda}(\infty)}
\newcommand{\vkinf}{v^{(k)}_{\lambda}(\infty)}
\newcommand{\C}{{\mathbb C}}
\newcommand{\Z}{{\mathbb Z}}  
\newcommand{\N}{{\mathbb N}}
\newcommand{\R}{{\mathbb R}}
\newcommand{\Q}{{\mathbb Q}}   
\newcommand{\F}{\mathcal F}
\newcommand{\Cc}{\mathcal C}
\newcommand{\Lc}{\mathcal L}
\newcommand{\mybox}{{\begin{picture}(15,15)%
\put(0,0){\line(1,0){15}}\put(0,0){\line(0,-1){15}}%
\put(15,0){\line(0,-1){15}}\put(0,-15){\line(1,0){15}}\end{picture}}}
\newcommand{\row}[2]{\multiput(0,#2)(15,0){#1}{\mybox}}
\newcommand{\bp}{\begin{picture}} \newcommand{\ep}{\end{picture}}
\newcommand{\rb}[1]{\raisebox{#1\unitlength}}
\renewcommand{\theequation}{\thesection.\arabic{equation}}
\def\theenumi{\roman{enumi}}
\def\labelenumi{(\theenumi)}

\title[Coinvariants and fusion product]
{Spaces of coinvariants and fusion product II.: 
{}$\slth$ character formulas in terms of Kostka polynomials}
\author{B. Feigin, M. Jimbo, R. Kedem, S. Loktev, T. Miwa}
\address{BF: Landau institute for Theoretical Physics, Chernogolovka,
142432, Russia}\email{feigin@feigin.mccme.ru}  
\address{MJ: Graduate School of Mathematical Sciences, The 
University of Tokyo, Tokyo 153-8914, Japan}\email{jimbomic@ms.u-tokyo.ac.jp}
\address{RK: Department of Mathematics, University of Illinois,
1409 W. Green St. Urbana, IL 61801, USA}
\email{rinat@math.uiuc.edu}
\address{SL:  Institute for Theoretical and Experimental Physics,
  B. Cheremushkinskaja, 25, Moscow 117259, Russia}
\address{
 Independent University of Moscow, B. Vlasievsky per, 11, 
Moscow 121002, Russia}\email{loktev@mccme.ru}
\address{TM: Division of Mathematics, Graduate School of Science, 
Kyoto University, Kyoto 606-8502
Japan}\email{tetsuji@kusm.kyoto-u.ac.jp}

\date{\today}

\title[Coinvariants and fusion product]
{Spaces of coinvariants and fusion product II.\\
{}$\slth$ character formulas in terms of Kostka polynomials}
\date{\today}
\begin{abstract}
  In this paper, we continue our study of the Hilbert polynomials of 
  coinvariants begun in our previous work \cite{FJKLM} (paper I).  
We describe the $\sln$ fusion products for symmetric tensor
representations following the method of \cite{FF}, and show that 
their Hilbert polynomials are $A_{n-1}$-supernomials. 
We identify the fusion product of arbitrary irreducible 
$\sln$-modules with the fusion product of 
their resctriction to $\mathfrak{sl}_{n-1}$. 
Then using the equivalence theorem from paper I and the results above
for $\mathfrak{sl}_3$ we give a fermionic formula for 
the Hilbert polynomials of a
  class of $\widehat{\mathfrak{sl}}_2$ coinvariants in terms of the
  level-restricted Kostka polynomials.  
The coinvariants under
  consideration are a generalization of the coinvariants 
 studied in \cite{FKLMM}. 
Our formula differs from the fermionic formula established in 
\cite{FKLMM} and implies the alternating sum formula 
conjectured in \cite{FL} for this case. 
\end{abstract}
\maketitle

\setcounter{section}{0}
\setcounter{equation}{0}
\section{Introduction}
Combinatorics and representation theory share a long history of mutual
interaction. Calculation of dimensions, multiplicities, etc., provides
interesting problems for combinatorics, and combinatorics gives a
natural language and technical machinery for solving them.  In many
cases, representation theory provides an interpretation of
combinatorial identities, or creates new identities.

The Kostka-Foulkes polynomial $K_{\lambda\mu}(q)$, the transition
matrix between Hall Littlewood and Schur functions, is one of the main players
in this interaction. In the original version, where $\lambda$ and $\mu$
are partitions, the Kostka number $K_{\lambda\mu}(1)$ is the cardinality
of the set of Young tableaux ${\mathcal T}(\lambda,\mu)$ of shape $\lambda$
and weight $\mu$. It gives the multiplicity of the irreducible
$\mathfrak{sl}_n$-module of highest weight $\lambda$ in the tensor product
of the $\mu_i$-th symmetric tensor representations $(1\leq i\leq N)$.

In combinatorics, introduction of $q$-polynomials such as $K_{\lambda\mu}(q)$
is equivalent to finding a good statistic on combinatorial objects.
In \cite{LS}, it was proved that the charge statistic on 
${\mathcal T}(\lambda,\mu)$ gives the Kostka polynomial.

In recent years, the interaction between combinatorics and representation
theory has been enhanced by integrable models in statistical mechanics
and quantum field theory. This is because physics provides the
equivalent of a charge statistic.

For the Kostka polynomial, \cite{KR} found a charge preserving
bijection between Young tableaux and rigged configurations.  The
latter appear naturally in the context of counting the Bethe Ansatz
states of the XXX Hamiltonian of one-dimensional quantum spin chain.
The meaning of the charge statistic for the Bethe states was later
clarified in \cite{KKMM}, where it was found that in the conformal
limit, the scaled partition function in the infinite size limit,
corresponding to the Bethe states, in many integrable models gives
characters in conformal field theory (CFT).  The bijection of
\cite{KR} led to new combinatorial identities which express the Kostka
polynomials as sums of products of $q$ binomial coefficients.  Such a
formula was called a fermionic formula in \cite{KKMM}.

Although a variety of fermionic formulas and related identities have
been obtained by combinatorial methods, the meaning from
representation theory was not given in general. In \cite{FJKLM} we
clarified some aspect of this connection between filtered spaces of
conformal blocks and the Kostka polynomial from representation theory
of affine Lie algebra $\widehat{\mathfrak g} ={\mathfrak
  g}\otimes{\mathbf C}[t,t^{-1}]\oplus{\mathbf C}K$. Namely, we
defined a natural filtration on the dual space of conformal blocks,
and showed that the Hilbert polynomial of the associated graded space
gives rise to the Kostka polynomial and its generalizations. Let us
recall this construction.

Let ${\mathcal Z}$ be a set of distinct complex numbers. 
To each $z_i$ we 
associate a left ideal $\gX_i\subset U({\mathfrak g})$. 
We fix a level $k$, and consider the finite-dimensional cyclic 
${\mathfrak g}$-module
\[
\pi^{(k)}(\gX_i)=\oplus_{\mu\in P^{(k)}_+}
\left(\pi_\mu^{\gX_i}\right)^*\otimes \pi_\mu.
\]
Here $\pi_\mu$ is the irreducible ${\mathfrak g}$-module 
with highest weight $\mu$, 
and $\left(\pi_\mu^{\gX_i}\right)^*$ is the dual 
space of $\gX_i$ invariants in $\pi_\mu$. 
The highest weight $\mu$ is restricted by the condition
$\left<\mu,\theta^\vee\right>\leq k$ where $\theta^\vee$ 
is the maximal coroot.

Let $L_i$ be a level $k$ integrable (in general, reducible) 
$\widehat{\mathfrak g}$-module 
localized at $z_i$, which has $\pi^{(k)}(\gX_i)$
as the degree zero part. 
Let us consider the space of conformal coinvariants
\bea
&&\left<L_0,L_1,\ldots,L_N\right>^{{\mathbf C}P^1}_{\infty,z_1,\ldots,z_N}\nn\\
&&=L_0\boxtimes L_1\boxtimes\cdots\boxtimes L_N/
{\mathfrak g}^{{\mathbf C}P^1\backslash\{\infty,z_1,\ldots,z_N\}}
(L_0\boxtimes L_1\boxtimes\cdots\boxtimes L_N),
\label{CFT}
\ena
where
$L_0=L^{(k)}_\lambda(\infty)$ is the level $k$ 
irreducible highest weight 
$\widehat{\mathfrak g}$-module with highest weight $\lambda$ 
`placed at $\infty$', whereas $L_i$ are placed at $z_i$
(see \cite{FJKLM} for the precise meaning of these terms). 
We fix a cyclic vector $v_i\in\pi^{(k)}(\gX_i)\subset L_i$
to be the canonical element corresponding to the inclusions
$\pi_\mu^{\gX_i}\rightarrow \pi_\mu$ $(\mu\in P^{(k)_+})$.
Set $X_i=S(\gX_i)+B_1$, 
where $B_1=\g\otimes t\C[t]U(\g[t])$, and 
$S$ is the antipode. Note that the left
ideal $S(X_i)$ annihilates $v_i$.
In the following, 
for a module $W$ over an algebra $A$ and a right ideal 
$Y\subset A$, we abbreviate the space of 
coinvariants $W/YW$ to $W/Y$. 

In \cite{FJKLM} we introduced the notion of 
the fusion right ideal 
$X_1\circledast\cdots\circledast X_N({\mathcal Z})$, 
and proved two canonical isomorphisms
\begin{equation}\label{ISO1}
L^{(k)}_\lambda(\infty)/X_1\circledast\cdots\circledast 
X_N({\mathcal Z})
\rightarrow
\left<L_0,L_1,\ldots,L_N\right>^{{\mathbf C}P^1}_{\infty,z_1,\ldots,z_N}
\end{equation}
and 
\begin{equation}\label{ISO2}
{\mathcal F}_{\mathcal Z}
\left(\pi^{(k)}(\gX_1),\ldots,\pi^{(k)}(\gX_N)\right)/
S\left(I^{(k)}_\lambda\right)
\rightarrow
\left<L_0,L_1,\ldots,L_N\right>^{{\mathbf C}P^1}_{\infty,z_1,\ldots,z_N}. 
\end{equation}
Here 
${\mathcal F}_{\mathcal Z}
\left(\pi^{(k)}(\gX_1),\ldots,\pi^{(k)}(\gX_N)\right)$
is the filtered tensor product (in the sense of \cite{FL}) of
$\pi^{(k)}(\gX_i)$
with the cyclic vector $v_1\otimes\cdots\otimes v_N$, and
$I^{(k)}_\lambda\subset U({\mathfrak g}[t])$ is a left ideal 
which
annihilates the highest weight vector 
$v_\lambda\in L^{(k)}_\lambda(\infty)$
(see \cite{FJKLM} for the precise definition).

The logical steps for the proof of these isomorphisms 
(\ref{ISO1}) and (\ref{ISO2}), given in \cite{FJKLM}, 
was as follows. 
A standard argument shows that there are canonical surjections.
The dimension of the space of conformal coinvariants 
(i.e., the right hand sides of (\ref{ISO1}) and (\ref{ISO2})) 
is given by the Verlinde rule. 
We showed that 
the dimension of the space in the left hand side of 
(\ref{ISO1}) obeys the same Verlinde rule. 
Then, we proved the equivalence 
theorem of filtered vector spaces
\begin{equation}\label{EQUIV}
L^{(k)}_\lambda(\infty)/X_1\circledast\cdots\circledast X_N({\mathcal Z})
\simeq
{\mathcal F}_{\mathcal Z}
\left(\pi^{(k)}(\gX_1),\ldots,\pi^{(k)}(\gX_N)\right)/
S\left(I^{(k)}_\lambda\right).
\end{equation}

In \cite{FJKLM}, we studied the case where $\g=\slt$ and 
the representations $L_i$ are irreducible. 
In that special case, 
we proved that the associated graded space of 
\eqref{ISO1}--\eqref{ISO2} is isomorphic to 
the space of coinvariants 
\bea
\pi^{(k)}(\gX_1)*\cdots*\pi^{(k)}(\gX_N)(\mathcal{Z})/
S\left(I^{(k)}_\lambda\right),
\label{FUS*}
\ena
where 	
$\pi^{(k)}(\gX_1)*\cdots*\pi^{(k)}(\gX_N)(\mathcal{Z})$ 
signifies the fusion product, i.e., 
the associated graded space of  
${\mathcal F}_{\mathcal Z}
\left(\pi^{(k)}(\gX_1),\ldots,\pi^{(k)}(\gX_N)\right)$. 
We also gave an explicit 
fermionic formula for the Hilbert polynomial
of \eqref{FUS*}.
In what follows we refer to the latter as 
the `character'.

The aim of the present article is 
to obtain similar results 
for a class of reducible representations $L_i$.
Let us describe the main points. 
Denote by $e,f,h$ the standard generators of $\g=\slt$.
In \cite{FJKLM}, we studied the space of coinvariants
\eqref{CFT}
by taking $\gX_i$ such that $\pi^{(k)}(\gX_i)$ is an 
irreducible ($l+1$)-dimensional module $\pi_l$ of $\g$. 
Here we choose $\gX_i$ of the form
\[
U(\g)e^{m_i+1}+U(\g)f\quad (1\leq m_i\leq k).
\]
In this case we have 
\begin{equation}\label{REDUC}
\pi^{(k)}(\gX_i)\simeq\oplus_{l=0}^{m_i}\pi_l,
\end{equation}
which is a reducible, cyclic module over $\g$.
We choose the cyclic vector to be the sum 
$\sum_{l=0}^{m_i}v_l$ 
of the lowest weight vectors $v_l\in\pi_l$.
In this setting, we compute the character of the 
associated graded space of the space of conformal 
coinvariants
\bea
\mathop{\rm gr}
\left<L_0,L_1,\ldots,L_N\right>^{{\mathbf C}P^1}_{\infty,z_1,\ldots,z_N}.
\label{gr}
\ena

For that purpose, we take the following approach. 
The direct sum (\ref{REDUC}) can be viewed as 
the $m_i$-th symmetric tensor representation
$\tilde{\pi}_{m_i}$ of $\sltr$ 
restricted to the subalgebra $\slt$.
We will show that the filtered tensor 
product appearing in (\ref{EQUIV}) for the reducible 
$\slt$-modules (\ref{REDUC}) is isomorphic to 
that for the irreducible 
$\sltr$-modules $\tilde{\pi}_{m_i}$ with the cyclic vector $v_{m_i}$.
This observation enables us to compute 
the right hand side of (\ref{EQUIV}).
As a result, 
we obtain a fermionic formula involving 
level-restricted Kostka polynomials for $\slt$
and $q$-binomial symbols. 
At $q=1$, this fermionic formula reduces to the known
formula for the dimension of the space of conformal blocks.
This establishes the isomorphism 
between \eqref{gr} and 
the space of coinvariants of the fusion product 
\eqref{FUS*}. 
The fermionic formula obtained above also 
settles a conjecture of \cite{FL}
in the special case $\g=\slt$
(see Theorem \ref{thm:3.2}). 
Our technical machinery from representation theory 
is not yet strong enough to prove the conjecture in general. 
We hope for a further progress in this direction.

In the case $m_1=\cdots=m_N=k$ 
with $z_1,\cdots,z_N$ being distinct,  
the left hand side of 
\eqref{EQUIV} reduces to the quotient of $L^{(k)}_\lambda(\infty)$ 
by the right ideal generated by 
\be
e\otimes t^i\prod_{j=1}^N(t-z_j),
\quad 
f\otimes t^i
\quad (i\ge 0).
\en
In \cite{FKLMM}, it was shown that its associated graded space 
coincides with the quotient of $L^{(k)}_\lambda(\infty)$ by the 
ideal of the above form with $z_1=\cdots=z_N=0$.  
Moreover a monomial basis of that space was constructed. 
The method of \cite{FKLMM} is based on 
the loop Heisenberg algebra contained in $\widehat{\sltr}$. 
As we mentioned already, the fermionic formula 
which naturally arises in this method 
is different from the one given in this paper. 
Though the two approaches are different, 
one could say that there is 
a common flavor in that both make use of the $\sltr$ structure.  
It would be an interesting problem 
to derive analogous results in the present setting as well. 

The plan of the paper is as follows. In Section \ref{sec:2},
we extend some result of \cite{FF} on the fusion product
in the case of $\sln$-modules. In particular, we obtain
the equality of the Hilbert polynomial of fusion product for
the symmetric tensor representations with the corresponding case of the
$A_{n-1}$ supernomials.
In Section \ref{sec:3}, we apply the result of Section \ref{sec:2}
for the $\sltr$ case in the computation of the
space of coinvariants in the $\slt$ fusion product.
Appendix is given for connecting filtered tensor products of
$\slnm$ and $\sln$.

\setcounter{section}{1}
\setcounter{equation}{0}

\section{Fusion product of modules over abelian Lie algebras}\label{sec:2}
In this section, we extend some of the results of \cite{FF} to the
$\sln$ setting for the case of symmetric tensor representations.
Namely we consider the fusion product of symmetric tensor
representations of $\sln$, and give an explicit characterization of
the annihilating ideals.  As in \cite{FF}, we derive short exact
sequences of these modules and obtain a decomposition of $\sln$ fusion
products into subquotients, each being isomorphic to an $\slnm$ fusion
product.  As a corollary, we identify their characters with the
$q$-supernomial coefficients.  We also give a recursive construction
of a basis of fusion products in terms of abelian currents.

\subsection{Symmetric tensors for $\sln$}\label{subsec:2.1}
We denote the generators of $\sln$ by $e_{ab},e_{ba}$ and
$h_{ab}=e_{aa}-e_{bb}$ ($1\le a<b\le n$), where
$e_{ab}=\bigl(\delta_{ia}\delta_{jb}\bigr)_{1\le i,j\le n}$ in the
defining representation.  We fix a
standard embedding $\mathfrak{sl}_2\subset\cdots\subset \slnm\subset
\sln$, regarding $\mathfrak{sl}_i$ as the subalgebra generated by
$e_{ab}$ ($n-i+1\le a\neq b\le n$).  In general, for a Lie algebra
$\g$ we set $\g[t]=\g\otimes\C[t]$ and write \be x[i]=x\otimes
t^i\quad (x \in \g, i\ge 0).  \en

Fix a set of distinct complex numbers 
\bea
\zz=(z_1,\cdots,z_N),\qquad z_i\neq z_j~~(i\neq j).
\label{zz} 
\ena
Consider the $k$-th symmetric tensor representation $W_k=S^k(\C^n)$  
with lowest weight vector $w_k$: $e_{ab}w_k=0$ for $a>b$. 
Our main concern is the fusion product \cite{FL}
\bea
W_{k_1}*\cdots *W_{k_N}(\zz), 
\label{fu1}
\ena
defined by choosing $w_{k_p}$ as the cyclic vector for each $W_{k_p}$.

A particularly simple feature of the symmetric tensors 
is that they can be described in terms of the 
abelian subalgebra 
\be
\ga=\C x_1\oplus\cdots\oplus\C x_{n-1}\subset\sln,\qquad
x_a=e_{n-a\,n}.
\en
Indeed, $W_k$ is generated by $w_k$ also as an $\ga$-module. 
Moreover, since $w_k$ is annihilated by  
all generators $e_{ab}$ ($a\neq b$) except for $x_1,\cdots,x_{n-1}$, 
the fusion product as $\ga[t]$-module 
is the same as the restriction to $\ga[t]$ of 
the fusion product \eqref{fu1} as $\sln[t]$-module.
For that reason 
we will henceforth forget the $\sln$-structure and 
focus attention to the $\ga$-structure. 

The universal enveloping algebra $U(\ga)$ is nothing but the 
polynomial ring $\C[x_1,\cdots,x_{n-1}]$. 
The annihilating ideal
\be
{\rm Ann}\,(W_k)=\{x\in\C[x_1,\cdots,x_{n-1}]\mid xw_k=0\}
\en
is generated by the elements $x_1^{\nu_1}\cdots x_{n-1}^{\nu_{n-1}}$ 
such that $\nu_1+\cdots+\nu_{n-1}=k+1$, $\nu_1,\cdots,\nu_{n-1}\ge 0$. 
One of our goals is to determine the annihilating ideal 
of the fusion product \eqref{fu1} 
in $U(\ga[t])=\C[x_a[i]\mid 1\le a\le n-1,i\in\Z_{\ge 0}]$. 
It turns out that $x_a[i]$ belongs to this ideal if $i$ is
sufficiently large (see \eqref{finite}), and hence 
we can work with polynomial rings with finitely many variables. 

In order to study the inductive structure of the fusion product,  
it is necessary to consider a family of modules 
slightly more general than \eqref{fu1}. 
In the next section, we begin by reformulating the setting to accommodate
this point. 

\subsection{Modules over an abelian Lie algebra}\label{subsec:2.2}

Let $\ga=\oplus_{a=1}^{n-1}\C x_a$ be an abelian Lie algebra. 
For $1\le m\le n$ and $k\ge 0$, set 
\bea
V^{(m)}_k=\C[x_1,\cdots,x_{n-1}]/J^{(m)}_k,
\label{mod1}
\ena
where $J^{(m)}_k$ is the ideal generated by the following elements:
\bea
J^{(m)}_k = \langle x_m,\cdots,x_{n-1},\quad
x_1^{\nu_1}\cdots x_{m-1}^{\nu_{m-1}}:
\nu_i\ge 0,\sum_{a=1}^{m-1}\nu_a=k+1\rangle.
\label{ide1}
\ena

Notice that 
$$
V_k^{(m)} \simeq {\underset{\nu_1 + ... + \nu_{m-1}\leq k}{\oplus}} 
\C x_1^{\nu_1} \cdots x_{m-1}^{\nu_{m-1}}
$$
so that 
$$
\dim V_k^{(m)} = \qbin{k+m-1}{m-1}.
$$

We regard \eqref{mod1} as a cyclic $\ga$-module with cyclic vector 
$v^{(m)}_k=1\bmod J^{(m)}_k$. 
It is a trivial module $\C\,1$ if $m=1$ or $k=0$.
Note also that $x_a=0$ on $V^{(m)}_k$ if $a\ge m$. 

Let $\zz$ be as in \eqref{zz}, and let
\be
\nb=(n_1,\cdots,n_N),
\quad
\kb=(k_1,\cdots,k_N)
\en
be $N$ tuples of non-negative integers with $1\le n_1,\cdots,n_N\le n$. 
We consider the filtered tensor product
\bea
\fil (V^{(n_1)}_{k_1},\cdots,V^{(n_N)}_{k_N})
=V^{(n_1)}_{k_1}(z_1)\otimes\cdots\otimes V^{(n_N)}_{k_N}(z_N)
\label{fil}
\ena
and the fusion product 
\bea
V_{\zz}(\nb,\kb)=V^{(n_1)}_{k_1}*\cdots * V^{(n_N)}_{k_N}(\zz). 
\label{fu2}
\ena
In \eqref{fil}, $V^{(n)}_{k}(z)$ denotes 
the evaluation module with parameter $z$, where 
$x[i]$ acts as $z^ix$ ($x\in\ga$). 
These are cyclic $\ga[t]$-modules generated by 
$v(\nb,\kb)=v^{(n_1)}_{k_1}\otimes\cdots\otimes v^{(n_N)}_{k_N}$.
We define a $\Z_{\ge 0}\times \Z_{\ge 0}^{n-1}$
grading on \eqref{fu2} by 
{\it degree} and {\it weight}, with the assignment
\be
&&
\deg x_a[i]=i,
\qquad
\wt x_a[i]=(0,\cdots,\overset{\scriptstyle{a{\rm -th}}}{1},\cdots,0),
\\
&&\deg v(\nb,\kb)=0,
\qquad
\wt v(\nb,\kb)=(0,\cdots,0).
\en
In general, for a graded vector space 
$W=\oplus_{d,m_1,\cdots,m_{n-1}}W_{d,m_1,\cdots,m_{n-1}}$,
we define its character by
\bea
&&{\rm ch}_{q,z_1,\cdots,z_{n-1}}W
=\sum q^d z_1^{m_1}\cdots z_{n-1}^{m_{n-1}}
\dim W_{d,m_1,\cdots,m_{n-1}}.
\label{chara}
\ena
By an inequality of characters 
${\rm ch}_{q,z_1,\cdots,z_{n-1}}W\le {\rm ch}_{q,z_1,\cdots,z_{n-1}}W'$
we mean inequalities for all the coefficients, i.e., 
$\dim W_{d,m_1,\cdots,m_{n-1}}\le \dim W'_{d,m_1,\cdots,m_{n-1}}$.

Given $\nb$, define $X_a$ to denote the set of representations on
which $x_a$ acts nontrivially, and $N_a$ to be the number of such
representations, that is,
\bea
&&X_a=\{p\mid 1\le p\le N,~~~n_p>a\}, 
\label{Xa}
\\
&&\Nb=(N_0,\cdots,N_{n-1}),
\quad N_a=|X_a|,
\label{Xa2}
\ena
so that $N=N_0\ge N_1\ge\cdots\ge N_n=0$. 
We note also that if $i\ge N_a$ then
$x_a[i]=0$ on the fusion product \eqref{fu2}.
This is because 
\begin{equation}\label{finite}
x_a\otimes t^{i-N_a}\prod_{p\in X_a}(t-z_p)
\end{equation}
acts as $0$ on the filtered tensor product \eqref{fil}.

\subsection{Annihilating ideal}\label{subsec:2.3}
Consider the polynomial ring 
$
R=\C[S_\Nb]
$
in indeterminates
\bea
S_\Nb=\{x_a[i]\mid 1\le a\le n-1,~0\le i\le N_a-1\}.
\ena
We denote by $R_{d,\mb}$ the homogeneous component of $R$ of 
degree $d$ and weight $\mb=(m_1,\cdots,m_{n-1})$.

The fusion product \eqref{fu2} has a presentation
as a quotient of $R$ 
\bea
V_{\zz}(\nb,\kb)=R/I_\zz(\nb,\kb)
\label{fu3}
\ena
by some homogeneous ideal $I_\zz(\nb,\kb)$. 
In this subsection, 
we will find a set of elements in this ideal which
are independent of $\zz$ (Proposition \ref{prop:2.1}). 
Later we show that they generate the whole 
ideal $I_\zz(\nb,\kb)$ (Theorem \ref{thm:2.1}). 

Set 
\bea
x_a(z)=\sum_{i=0}^{N_a-1}x_a[i]z^{N_a-1-i}.
\label{xaz}
\ena
For $\nub=(\nu_1,\cdots,\nu_{n-1})\in \Z^{n-1}_{\ge 0}$, 
define
\bea
\mu(\nub,\kb)=\sum_{p=1}^N
\Bigl(\sum_{a=1}^{n_p-1}\nu_a-k_p\Bigr)_+
\label{mu}
\ena
where $u_+=\max(u,0)$. 
Consider the ideal $J(\nb,\kb)$ of $R$ generated as:
\begin{equation}
J(\nb,\kb) = \langle
\mbox{coefficients of $z^j$ in 
$\prod_{a=1}^{n-1}x_a(z)^{\nu_a}$},\ \forall j<\mu(\nub,\kb),  \nub\in \Z^{n-1}_{\ge 0}
\rangle. 
\label{J}
\end{equation}

\begin{prop}\label{prop:2.1}
We have 
\be
J(\nb,\kb)\subset I_{\zz}(\nb,\kb).
\en
That is, the following relations hold in $R/I_{\zz}(\nb,\kb)$:
\bea
\ord_z\prod_{a=1}^{n-1}x_a(z)^{\nu_a}\ge \mu(\nub,\kb).
\label{ordz}
\ena
\end{prop}
Above, for a polynomial $\varphi(z)=\sum_{j\ge 0}c_jz^j$ 
we write $\ord_z\varphi\ge j_0$ if $c_j=0$ for $j<j_0$.

To prove Proposition \ref{prop:2.1}, we prepare a few Lemmas.
For each $0\le a\le n-1$, set 
\be
\tilde{x}_a(z)
&=&\sum_{i=0}^{N_a-1}\tilde{x}_a[i]z^{N_a-1-i},
\\
\tilde{x}_a[i]
&=&
x_a[i]-x_a[i-1]\sigma_{a,1}+x_a[i-2]\sigma_{a,2}-\cdots
\\
&=&\sum_{s=0}^{i}(-1)^s x_a[i-s]\sigma_{a,s},
\en
where $\sigma_{a,s}$ denotes the $s$-th elementary symmetric polynomial
in the variables $\{z_p\}_{p\in X_a}$. 
Let $x_a^{(p)}$ stand for the action of $x_a$ on the 
$p$-th tensor component. 
\begin{lem}\label{lem:2.1}
There exists a constant $c^{(p)}_a\neq 0$ such that 
on \eqref{fil} we have
\bea
\tilde{x}_a(z_p)=c^{(p)}_ax^{(p)}_a
\qquad (1\le a\le n_p-1).
\label{xtil}
\ena
\end{lem}
\begin{proof}
Let $v$ be a vector of \eqref{fil}. Then we have
\be
\sum_{i=0}^\infty z^{-i}\tilde{x}_a[i]v
&=&
\prod_{p\in X_a}(1-z_p/z)\sum_{i=0}^{\infty} z^{-i}x_a[i]v
\\
&=&
\prod_{p\in X_a}(1-z_p/z)\sum_{p\in X_a}\frac{1}{1-z_p/z}x_a^{(p)}v
\\
&=&\sum_{p\in X_a}\prod_{p'(\neq p)\atop p'\in X_a}
(1-z_{p'}/z)\cdot x_a^{(p)}v.
\en
In the second line, we used $x^{(p)}_av=0$ for $p\not\in X_a$. 
We have shown that $\tilde{x}_a[i]v=0$ for $i\ge N_a$. 
Multiplying $z^{N_a-1}$ and setting $z=z_p$, we obtain
\eqref{xtil} with 
$c^{(p)}_a=\prod_{p'(\neq p)\atop p'\in X_a}(z_p-z_{p'})$.
\end{proof}

Let us present the filtered tensor product as 
a quotient of $R$ 
by some ideal $\tilde{I}_\zz(\nb,\kb)$,
\be
\fil (V^{(n_1)}_{k_1},\cdots,V^{(n_N)}_{k_N})
=R/\tilde{I}_\zz(\nb,\kb). 
\en
Recall from \eqref{ide1} that
\[
\prod_{a=1}^{n_p-1}(x^{(p)}_a)^{\nu_a}\in J^{(n_p)}_{k_p}
\quad \hbox{ if }\quad \sum_{a=1}^{n_p-1}\nu_p>k_p.
\]
We will now interpret this condition in terms of $x_a[i]$ by using the current
$x_a(z)$.

Denote by $R_d=\oplus_{\mb}R_{d,\mb}$ 
the homogeneous component of degree $d$, 
and set $R_{\le d}=R_0+\cdots+R_d$. 
\begin{lem}\label{lem:2.2}
(i)The ideal $I_{\zz}(\nb,\kb)$ consists of 
homogeneous elements $P\in R_d$ ($d\ge 0$), such that 
\be
P\equiv P' \bmod R_{\le d-1}
\qquad 
\mbox{for some $P'\in \tilde{I}_\zz(\nb,\kb)\cap R_{\le d}$.}
\en
(ii)For each $p$, the elements
\be
\prod_{a=1}^{n_p-1}\tilde{x}_a(z_p)^{\nu_p},
\qquad (\nub\in\Z^{n-1}_{\ge 0},\sum_{a=1}^{n_p-1}\nu_a>k_p)
\en
belong to the ideal $\tilde{I}_\zz(\nb,\kb)$.
\end{lem}
\begin{proof}
Assertion (i) is clear from the definition of the fusion product as
the adjoint graded space of the filtered tensor product. 
Assertion (ii) is an immediate consequence of Lemma \ref{lem:2.1} and
\ref{ide1}.
\end{proof}
\medskip

\noindent{\it Proof of Proposition \ref{prop:2.1}.}
\quad 
Let us prove the relation \eqref{ordz}. 
Fix $\nub$, and consider the polynomial
\be
f(z)=\prod_{a=1}^{n-1}\tilde{x}_a(z)^{\nu_a}\quad \in R[z].
\en
{}From Lemma \ref{lem:2.2} (ii), we have for each $p=1,\cdots,N$
\bea
\left.\frac{d^jf(z)}{dz^j}\right|_{z=z_p}\in \tilde{I}_{\zz}(\nb,\kb)
\qquad (0\le j\le \mu_p-1),
\label{aaa}
\ena
where 
\be
\mu_p=\bigl(\sum_{a=1}^{n_p-1}\nu_a-k_p\bigr)_+.
\en

Using \eqref{aaa} and the fact that $z_1,\cdots,z_N$ are distinct, 
we can find polynomials $g(z),h(z)\in R[z]$ such that
\bea
&&f(z)=\prod_{p=1}^N(z-z_p)^{\mu_p}\cdot g(z)+h(z),
\label{fgh}
\\
&&h(z)\in \tilde{I}_{\zz}(\nb,\kb)[z],
\quad \deg h<\deg g. 
\nn
\ena
Note that they have the form
\be
f(z)=\sum_{j=0}^d a_jz^{d-j},
\quad
g(z)=\sum_{j=0}^{d'} b_jz^{d'-j},
\quad 
h(z)=\sum_{j=d-d'+1}^{d} c_{j} z^{d-j},
\en
where $a_j,b_j,c_j\in R_{\le j}$. 

Let $\varphi_t$ be the automorphism of $R[t,t^{-1}]$ 
defined by $\varphi_t(b)=t^{-i}b$ for $b\in R_i$.
Because of the above form, 
$t^d\varphi_t(f)(z/t)$, 
$t^{d'}\varphi_t(g)(z/t)$ and $t^d\varphi_t(h)(z/t)$ 
have a well-defined limit as $t\rightarrow 0$.
Denote the limit by $f^*(z),g^*(z), h^*(z)$, respectively. 
Clearly we have
\be
f^*(z)=\prod_{a=1}^{n-1}x_a(z)^{\nu_a},
\qquad 
h^*(z)\in I_{\zz}(\nb,\kb).
\en
{}From \eqref{fgh} we obtain
\be
f^*(z)=z^{\sum_{p=1}^N\mu_p}g^*(z)+h^*(z).
\en
The proof is over.
\qed

\subsection{An isomorphism}\label{subsec:2.4}
Consider now the module 
\bea
W(\nb,\kb)=
W\!\!\begin{pmatrix}n_1,&\cdots,&n_N\\ k_1,&\cdots,&k_N\\ 
\end{pmatrix} 
=R/J(\nb,\kb).
\label{W}
\ena
The following relations are immediate from the definition.
\bea
&&
W\!\!
\begin{pmatrix} 
\cdots &n_{p-1} &1    &n_{p+1}&\cdots  \\
\cdots &k_{p-1} &k_p  &k_{p+1}&\cdots  \\
\end{pmatrix}
=
W\!\!
\begin{pmatrix}
\cdots &n_{p-1} &n_{p+1}&\cdots  \\
\cdots &k_{p-1} &k_{p+1}&\cdots  \\
\end{pmatrix},
\label{n=1}
\\
&&W\!\!
\begin{pmatrix}
\cdots &n_p     & \cdots  &n_{p'}&\cdots  \\
\cdots &k_p     & \cdots  &k_{p'}&\cdots  \\
\end{pmatrix}
=
W\!\!
\begin{pmatrix}
\cdots &n_{p'}& \cdots  &n_p   &\cdots  \\
\cdots &k_{p'}& \cdots  &k_p   &\cdots  \\
\end{pmatrix}.
\label{Wsym}
\ena
We have also a canonical isomorphism
\bea
&&
W\!\!
\begin{pmatrix}
\cdots &n_{p-1} &n_p&n_{p+1}&\cdots  \\
\cdots &k_{p-1} &0  &k_{p+1}&\cdots  \\
\end{pmatrix}
\simeq
W\!\!
\begin{pmatrix}
\cdots &n_{p-1} &n_{p+1}&\cdots  \\
\cdots &k_{p-1} &k_{p+1}&\cdots  \\
\end{pmatrix}.
\label{k=0}  
\ena
To see this note that if $k_p=0$, then for all $a<n_p$, 
the generator $x_a[N_a-1]$
belongs to the ideal $J(\nb,\kb)$
(take $\nu_b=\delta_{ab}$).
Therefore, in the presentation \eqref{W}, 
we can replace $R$ by dropping all such
$x_a[N_a-1]$.
The result is the presentation in the right hand side.
In what follows we identify modules which are related 
by \eqref{k=0}.

Our first goal is the following result.
\begin{thm}\label{thm:2.1}
We have an isomorphism 
\bea
W(\nb,\kb) \overset{\sim}{\longrightarrow} V_{\zz}(\nb,\kb)
\label{surj1}
\ena
given by the canonical surjection.
In particular, the fusion product is independent of the choice of 
$\zz$.
\end{thm}
We prove Theorem \ref{thm:2.1} in subsection 
\ref{subsec:2.8}.

\subsection{Subquotient modules}\label{subsec:2.5}
In order to determine the structure of 
$W(\nb,\kb)$, we study its subquotients. 
If $k_p=0$ for all $p$ such that $n_p\ge 2$, then 
$W(\nb,\kb)\simeq\C$.  
Otherwise, we can find a number $p$ satisfying the following conditions.
\bea
&&n_p\ge 2,~~k_p>0,
\quad
\mbox{ $k_s\le k_p$ for $n_s\ge n_p$}.
\label{order3}
\ena
In what follows we fix $p$ subject to the conditions \eqref{order3} and set 
\bea
a=n_p-1, 
\quad 1\le a\le n-1.
\label{order4}
\ena

Define 
\be
&&\nb'=(n_1,\cdots,n_p-1,\cdots,n_N), 
\\
&&\Nb'=(N_0,\cdots,N_a-1,\cdots,N_{n-1}), 
\\
&&\kb''=(k_1,\cdots,k_p-1,\cdots,k_N). 
\en
We set $S=S_\Nb$, 
$S'=S_{\Nb'}$, $S''=S_{\Nb}$, 
$J=J(\nb,\kb)$, 
$J'=J(\nb',\kb)$, 
$J''=J(\nb,\kb'')$ and 
\be
&&
W'=W\!\!\begin{pmatrix}
n_1&\cdots&n_p-1&\cdots&n_N\\
k_1&\cdots&k_p&\cdots&k_N\\
\end{pmatrix}
=\C[S']/J', 
\\
&&
W=
W\!\!\begin{pmatrix}
n_1&\cdots &n_p&\cdots&n_N\\
k_1&\cdots &k_p&\cdots&k_N\\
\end{pmatrix}
=\C[S]/J, 
\\
&&
W''=
W\!\!\begin{pmatrix}
n_1&\cdots &n_p&\cdots&n_N\\
k_1&\cdots &k_p-1&\cdots&k_N\\
\end{pmatrix}
=\C[S'']/J''.
\en
Our second goal is the following theorem. 
\begin{thm}\label{thm:2.2}
Notation being as above, 
there exists an exact sequence 
\bea
&&
0\longrightarrow W'\overset{\iota}{\longrightarrow} W
\overset{\psi}{\longrightarrow}
W''\longrightarrow 0,
\label{exact1}
\ena
where the maps shift the grading as 
\be
&&\iota\left(W'_{d-m_a,m_1,
\cdots,m_a,\cdots,m_{n-1}}\right)\subset 
W_{d,m_1,\cdots,m_a,\cdots,m_{n-1}},\\
&&\psi\left(W_{d,m_1,\cdots,m_a,\cdots,m_{n-1}}\right)\subset 
W''_{d,m_1,\cdots,m_a-1,\cdots,m_{n-1}}.
\en
\end{thm}
The maps $\iota$ and $\psi$ are described as follows. 

Let $\tilde{\iota}$ be the map 
\be
\tilde{\iota}:\C[S']\rightarrow \C[S],
\quad 
x_b[i]\mapsto x_b[i+\delta_{ab}]. 
\en
In subsection \ref{subsec:2.6} we show 
that $\tilde{\iota}(J')\subset J$ (Proposition \ref{prop:2.2}).  
The map $\iota$ is the induced map 
\bea
\iota~:~W'=\C[S']/J' \longrightarrow W=\C[S]/J.
\label{iota}
\ena

Set $\overline{W}'=\iota(W')\subset W$.  
In other words, $\overline{W}'$ is the subspace spanned by all monomials in 
$x_b[i]$ ($b\neq a$, $0\le i\le N_b-1$) and 
$x_a[1],\cdots,x_a[N_a-1]$. 
The quotient $W/\overline{W}'$ is spanned by 
monomials which are divisible by $x_a[0]$. 

In subsection \ref{subsec:2.7}, we show that there is a well-defined map 
\bea
\phi~:~W''\longrightarrow W/\overline{W}'
\label{phi}
\ena
given by
\be
x_a[0]^i b~~\mapsto~~\frac{1}{i+1}x_a[0]^{i+1} b,
\en
where $b$ runs over monomials not divisible by $x_a[0]$. 
We show moreover that \eqref{phi} is an isomorphism. 
The map $\psi$ is the composition $\phi^{-1}\circ\pi$, 
where $\pi:W\rightarrow W/\overline{W}'$ signifies the canonical surjection.

Theorem \ref{thm:2.2} will be proved in subsection \ref{subsec:2.8} 
simultaneously with Theorem \ref{thm:2.1}.

\subsection{The map $\iota$}\label{subsec:2.6}
Here we show that the map \eqref{iota} is
well-defined.

\begin{prop}\label{prop:2.2}
We have $\tilde{\iota}(J')\subset J$.  
\end{prop}
\begin{proof}
Set 
\be
x_b'(z)=\sum_{i=0}^{N_b'-1}x_b[i]z^{N_b'-1-i},
\quad
N_b'=N_b-\delta_{ab}.
\en
Then 
$\tilde{\iota}(x_b'(z))=x_b(z)-\delta_{ab}x_a[0]z^{N_a-1}$, 
so that we have
\bea
&&\tilde{\iota}\left(\prod_{b=1}^{n-1}x_b'(z)^{\nu_b}\right)
\label{aux}\\
&&\quad =
\sum_{\mu=0}^{\nu_a}(-x_a[0])^\mu
\prod_{b(\neq a)}x_b(z)^{\nu_b}\cdot
x_a(z)^{\nu_a-\mu}z^{(N_a-1)\mu}.
\nn
\ena
Let $C(j)$ denote the coefficient of $z^j$ in 
$\prod_{b=1}^{n-1}x'_b(z)^{\nu_b}$. 
We are to show that $\tilde{\iota}\bigl(C(j)\bigr)$
belongs to $J$ if $j<d_1$, where 
\be
d_1=\sum_{s=1}^{N}
\Bigl(\sum_{b=1}^{n'_s-1}\nu_b-k_s \Bigr)_+,
\quad n_s'=n_s-\delta_{sp}.
\en

Let $C(\mu,j)$ denote the coefficient of $z^j$ in 
the $\mu$-th term in the right hand side of \eqref{aux}.
It belongs to $J$ if $j<d_2$, where
(see \eqref{Xa} and \eqref{Xa2} for $X_a$)
\be
d_2=\sum_{s\in X_a}
\Bigl(\sum_{b=1}^{n_s-1}\nu_b-\mu-k_s\Bigr)_+
+\sum_{s\not\in X_a}
\Bigl(\sum_{b=1}^{n_s-1}\nu_b-k_s\Bigr)_+
+(N_a-1)\mu.
\en
Noting that $x_++y\ge (x+y)_+$ holds for 
$y\ge 0$, $a=n_p-1$, 
and that $\mu\leq\nu_a=\nu_{n_p-1}$ implies
\[
\sum_{b=1}^{n_p-1}\nu_b-\mu-k_p
\ge
\sum_{b=1}^{n_p-2}\nu_b-k_p=
\sum_{b=1}^{n'_p-1}\nu_b-k_p,
\]
we find $d_2\ge d_1$. The assertion follows from this.
\end{proof}

\subsection{Dual space}\label{subsec:2.7}
In order to construct the map $\phi$ in \eqref{phi},  
we pass to the dual space $W(\nb,\kb)^*$ and use 
a realization for the latter in terms of polynomials. 

Recall $R=\C[S_{\bf N}]$, 
with $S_{\bf N}=\{x_b[i]\mid 0\le i\le N_b-1,
1\le b\le n-1\}$. 
Consider the restricted dual space
\be
R^*=\bigoplus_{d,\mb}(R_{d,\mb})^*.
\en
Each homogeneous component $(R_{d,\mb})^*$
is realized as a subspace of the ring of 
polynomials $\C[T_\mb]$ in indeterminates 
\be
T_\mb=\{y_{bj}\mid 1\le j\le m_b,~~ 1\le b\le n-1\}, 
\en
where $\mb=(m_1,\cdots,m_{n-1})$. 
In this space
we count the total degree $d^*$ by assigning $\deg y_{bj}=1$. 
Let $\mathcal{F}_{d^*,\mb,\Nb}$ denote the subspace of 
 $\C[T_\mb]$ consisting of all polynomials 
$f$ with the properties 
\begin{enumerate}
\item $f$ is symmetric in 
each group of variables $(y_{b1},\cdots,y_{bm_b})$, 
\item $f$ has total degree $d^*$, 
and $\deg_{y_{bj}}f<N_b$ for each $b,j$.
\end{enumerate}
Set 
\be
d(\Nb,\mb)=\sum_{b=1}^{n-1}(N_b-1)m_b. 
\en
We have an isomorphism 
\bea
R_{d,m})^*
\overset{\sim}{\longrightarrow}
\mathcal{F}_{d^*,\mb,\Nb},
\quad d+d^*=d(\Nb,\mb)
\label{ST}
\ena
given by
\be
\bigoplus_{d=0}^{d(\Nb,\mb)}
\C[S_{\bf N}]_{d,\mb}^*
~~
\ni 
\theta\mapsto 
\theta\Bigl(
\prod_{1\le b\le n-1\atop 1\le j\le m_b}
x_b\bigl(y_{bj}\bigr)
\Bigr) 
~~
\in
\bigoplus_{d^*=0}^{d(\Nb,\mb)}
\mathcal{F}_{d^*,\mb,\Nb}.
\en
With respect to this pairing, monomials 
\be
\prod_{a=1}^{n-1}x_a[i^{(a)}_1]\cdots x_a[i^{(a)}_{m_a}]
\quad
(0\le i^{(a)}_1\le \cdots\le i^{(a)}_{m_a}\le N_a-1)
\en
and 
\be
\prod_{a=1}^{n-1}
\mathop{\rm Sym}\bigl(y_{a1}^{N_a-1-i^{(a)}_1}\cdots
y_{am_a}^{N_a-1-i^{(a)}_{m_a}}\bigr)
\en
are dual bases to each other, 
where the symbol 
$\mathop{\rm Sym}$
stands for the normalized sum over all permutations of 
$i^{(a)}_1,\cdots,i^{(a)}_{m_a}$. 

We have introduced the quotient space
$W(\nb,\kb)$ of $R=\C[S_\Nb]$ by the homogeneous ideal 
$J=J(\nb,\kb)$ \eqref{J}.
The dual space $\bigl(W(\nb,\kb)_{d,\mb}\bigr)^*$ 
is the orthogonal complement to $J$ under the 
coupling given above. 
Explicitly it is described as follows. 
Given $\nub$ such that $\nu_b\le m_b$ ($1\le b\le n-1$), 
we consider the specialization of variables 
\bea
\begin{matrix}
y_{11}=&\cdots&=y_{1\nu_1}=z\\
       &\cdots&           \\    
y_{n-1\,1}=&\cdots&=y_{n-1\,\nu_{n-1}}=z\\
\end{matrix}
\label{special}
\ena
and denote by $\Phi(f,\nu_1,\cdots,\nu_{n-1};z)$ 
the polynomial obtained from $f$ by \eqref{special}. 

\begin{lem}\label{lem:2.4}
The dual space 
$\bigl(W(\nb,\kb)_{d,\mb}\bigr)^*$
is isomorphic to the subspace of 
$\mathcal{F}_{d^*,\mb,\Nb}$ 
consisting of polynomials $f$ satisfying $($see \eqref{mu} for 
$\mu(\nub,\kb)$$)$
\bea
\ord_z\Phi(f,\nu_1,\cdots,\nu_{n-1};z)\ge 
\mu(\nub,\kb)
\label{ortho}
\ena
for all $\nub\in\Z_{\ge 0}^{n-1}$ such that
$\nu_b\le m_b$ $(1\le b\le n-1)$. 
Namely, $f\in \mathcal{F}_{d^*,\mb,\Nb}$ 
is orthogonal to $J$ if and only if it satisfies the condition \eqref{ortho}.
\end{lem}

We need also the orthogonal complement 
to the subspace ${\rm Im}\,\tilde{\iota}$.
An element $f\in \mathcal{F}_{d^*,\mb,\Nb}$ is orthogonal 
to ${\rm Im}\,\tilde{\iota}$ if and only if 
\bea
\label{NOMO}
&&\hbox{it does not contain monomials of the form}
\prod_{1\le b\le n-1\atop 1\le j\le m_b}y_{bj}^{i_{bj}}\\
&&\quad \mbox{ with }i_{a1},\cdots,i_{am_a}<N_a-1.
\nn
\ena
Such an $f$  can be written uniquely in the form 
\bea
f=\sum_{j=1}^{m_a}y_{aj}^{N_a-1}g^{(j)}
\label{fg}
\ena
with some polynomial 
$g\in\C[T_{m_1,\cdots,m_a-1,\cdots,m_{n-1}}]$, where 
\be
g^{(j)}=
g(y_{11},\cdots;\cdots,\overset{\frown}{y_{aj}},\cdots;
\cdots,y_{n-1\,m_{n-1}}).
\en
In the right hand side the variable $y_{aj}$ is omitted.
Let 
$\mathcal{F}'_{d^*,\mb,\Nb}$ 
be the subspace of polynomials satisfying \eqref{NOMO}. 
Setting $\tilde{\phi}^*(f)=g$, we obtain an injective map
\be
\tilde{\phi}^*~:~
\mathcal{F}'_{d^*,m_1,\cdots,m_a,\cdots,m_{n-1},\Nb} 
\longrightarrow
\mathcal{F}_{d^*-N_a+1,m_1,\cdots,m_a-1,\cdots,m_{n-1},\Nb}. 
\en
Note that the dual space 
$(W/\overline W')^*_{d^*,{\mathbf m}}$ is isomorphic to the subspace of 
${\mathcal F}'_{d^*,{\mathbf m},{\mathbf N}}$ which is orthogonal to $J$.

\begin{prop}\label{prop:2.3}
If $f\in \mathcal{F}'_{d^*,\mb,\Nb}$ is orthogonal to $J$, then 
$g=\tilde{\phi}^*(f)$ is orthogonal to $J''$. 
Hence $\tilde{\phi}^*$ induces an injection
\be
\phi^*~:~
(W/\overline W')^*_{d^*,m_1,\cdots,m_{a},\cdots,m_{n-1}}
\longrightarrow 
(W'')^*_{d^*-N_a+1,m_1,\cdots,m_a-1,\cdots,m_{n-1}}. 
\en
\end{prop}
\begin{proof}
By Lemma \ref{lem:2.4}, $f$ satisfies the condition \eqref{ortho}.
We are to show that 
\bea
\ord_z\Phi(g,\nu_1,\cdots,\nu_{n-1};z)
\ge 
\sum_{s=1}^N
\bigl(
\sum_{1\le b\le n_s-1}\nu_b-k_s+\delta_{sp}\bigr)_+
\label{aim}
\ena
holds for any $\nub$ with $\nu_b\le m_b-\delta_{ab}$ ($1\le b\le n-1$). 

Specializing \eqref{fg}, we obtain 
\bea
&&\Phi(f,\nu_1,\cdots,\nu_a,\cdots,\nu_{n-1};z)
\label{induc}\\
&&=\nu_a z^{N_a-1}
\Phi(g,\nu_1,\cdots,\nu_a-1,\cdots,\nu_{n-1};z)
\nn\\
&&+\sum_{\nu_a+1\le j\le m_a}y_{aj}^{N_a-1}
\Phi(g^{(j)},\nu_1,\cdots,\nu_a,\cdots,\nu_{n-1};z). 
\nn
\ena
We distinguish two cases. 
\medskip

\noindent{{\bf Case $\sum_{1\le b\le a}\nu_b \le k_p-1$}:}
\quad 
In this case, the term with $s=p$ 
in the right hand side of \eqref{aim}
does not contribute. 
If $\nu_a=0$, then \eqref{aim} follows from 
\eqref{ortho}, 
\eqref{induc} and the 
uniqueness of the representation \eqref{fg}.
The general case follows by ascending induction on $\nu_a$. 
Namely, we apply the induction hypothesis 
\eqref{aim}
for the first term
in the right hand side of 
\eqref{induc}, 
and conclude 
\eqref{aim}
for the second term.

\noindent{{\bf Case $\sum_{1\le b\le a}\nu_b \ge k_p$}:}
\quad  
We use (\ref{induc}) by changing $\nu_a$ to $\nu_a+1$:
\begin{eqnarray*}
&&\Phi(f,\nu_1,\cdots,\nu_a+1,\cdots,\nu_{n-1};z)
\\
&&=(\nu_a+1)z^{N_a-1}
\Phi(g,\nu_1,\cdots,\nu_a,\cdots,\nu_{n-1};z)
\nn\\
&&+\sum_{\nu_a+2\le j\le m_a}y_{aj}^{N_a-1}
\Phi(g^{(j)},\nu_1,\cdots,\nu_a+1,\cdots,\nu_{n-1};z). 
\nn
\end{eqnarray*}
Consider first the case $\nu_a=m_a-1$. 
The second term in the right hand side is then void. 
We have 
\bea
&&\ord_z\Phi(g,\nu_1,\cdots,\nu_a,\cdots,\nu_{n-1};z)
\nn\\
&&\ge
\sum_{s\not\in X_a}
\bigl(\sum_{1\le b\le n_s-1}\nu_b-k_s\bigr)_+
+
\sum_{s\in X_a}
\bigl(\sum_{1\le b\le n_s-1}\nu_b+1-k_s\bigr)_+
-(N_a-1).
\label{xxxx}
\ena
Since $a=n_p-1$, $s\in X_a$ implies $n_p\le n_s$, and hence 
\be
\sum_{1\le b\le n_s-1}\nu_b
\ge
\sum_{1\le b\le a}\nu_b
\ge
k_p
\ge k_s.
\en
In the last step we used the assumption \eqref{order3}. 
Therefore we find that the right hand side of \eqref{xxxx}
is not less than
\be
\sum_{s=1}^N
\bigl(
\sum_{1\le b\le n_s-1}\nu_b-k_s+\delta_{sp}\bigr)_+.
\nn
\en
The general case follows by descending induction on $\nu_a$. 
Namely, we apply the induction hypothesis 
\eqref{aim}
for the second term
in the right hand side of 
\eqref{induc}, 
and conclude 
\eqref{aim}
for the first term.
\end{proof}

Let 
\be
\tilde{\phi}~:~ {\bf C}[S_{\bf N}]_{d,m_1,\ldots,m_a-1,\ldots,m_{n-1}}
\longrightarrow 
({\bf C}[S_{\bf N}]/{\rm Im}\,\tilde{\iota})_{d,{\bf m}}
\en
be the map dual to $\tilde{\phi}^*$. 
Decomposing ${\bf C}[S_{\bf N}]_{d,m_1,\ldots,m_a-1,\ldots,m_{n-1}}$ as 
\[
\oplus_{i=0}^{m_a-1}x_a[0]^i
({\rm Im}\,\tilde{\iota})_{d,m_1,\ldots,m_a-1-i,\ldots,m_{n-1}}, 
\]
we can explicitly write $\tilde{\phi}$ as
\[
\tilde{\phi}:x_a[0]^i({\rm Im}\,\iota)_{d,m_1,\ldots,m_a-1-i,\ldots,m_{n-1}}
\ni b\mapsto \frac{x_a[0]}{i+1}b\in
({\bf C}[S_{\bf N}]/{\rm Im}\,\tilde{\iota})_{d,{\bf m}}.
\]
We have shown in Proposition \ref{prop:2.3} that if 
$b\in{\bf C}[S_{\bf N}]_{d,m_1,\ldots,m_a-1,\ldots,m_{n-1}}$ belongs to
$J''$ then the image $\tilde{\phi}(b)$ belongs to $J$. 
Therefore, $\tilde{\phi}$ induces the 
surjection $\phi$ of \eqref{phi}.

\subsection{Proof of Theorems \ref{thm:2.1},\ref{thm:2.2}}\label{subsec:2.8}
Summarizing the results of subsections \ref{subsec:2.6} and \ref{subsec:2.7}, 
we obtain the diagram 
\be
&&
\begin{CD}
\phantom{0}@.W'@.@.W''@.
\\
@.@V{\iota}VV @.\phantom{W}@VV{\phi}V @.@.
\\
0@>>>\overline{W}' @>>> W @>\pi>> W/\overline{W}' @>>>0
\end{CD}
\en
where the vertical arrows are surjective. 
In particular we have
\bea
\dim W\le \dim W'+\dim W''.
\label{est1}
\ena

We are now in a position to finish the proof 
of Theorems \ref{thm:2.1}, \ref{thm:2.2}.
Set 
\be
&&\tilde{d}(\nb,\kb)=\dim W(\nb,\kb),
\\
&&
d(\nb,\kb)=\dim V_{\zz}(\nb,\kb).
\en
{}From \eqref{surj1} we have 
$\tilde{d}(\nb,\kb)\ge d(\nb,\kb)$.

On the other hand, under the assumption \eqref{order3}, 
we have from \eqref{est1} the inequality 
\bea
&&
\label{ineq}
\\
&&
\tilde{d}
\begin{pmatrix}
\cdots&a+1&\cdots \\
\cdots&k_p&\cdots \\
\end{pmatrix}
\le
\tilde{d}
\begin{pmatrix}
\cdots&a&\cdots \\
\cdots&k_p&\cdots \\
\end{pmatrix}
+
\tilde{d}
\begin{pmatrix}
\cdots &a+1&\cdots\\
\cdots &k_p-1&\cdots\\
\end{pmatrix}.
\nn
\ena
We have also the relations  
\be
&&
\tilde{d}\!\!
\begin{pmatrix}
1,&\cdots,&1\\ 
k_1,&\cdots,&k_N\\ 
\end{pmatrix}=1,
\\
&&
\tilde{d}\!\!
\begin{pmatrix}
\cdots&n_p&\cdots\\
\cdots&0&\cdots\\
\end{pmatrix}
=
\tilde{d}\!\!
\begin{pmatrix}
\cdots&n_p-1&\cdots\\
\cdots&0&\cdots\\
\end{pmatrix},
\\
&&
\tilde{d}
\begin{pmatrix}
\cdots &n_p     & \cdots  &n_{p'}&\cdots  \\
\cdots &k_p     & \cdots  &k_{p'}&\cdots  \\
\end{pmatrix}
=\tilde{d}
\begin{pmatrix}
\cdots &n_{p'}& \cdots  &n_p   &\cdots  \\
\cdots &k_{p'}& \cdots  &k_p   &\cdots  \\
\end{pmatrix}.
\en
The above relations are also true for $d(\nb,\kb)$,
and moreover the equality holds in \eqref{ineq}, due to the identity
for binomial coefficients. 
These equations determine $d(\nb,\kb)$ uniquely.  
Therefore by induction we have for all $\nb,\kb$ that 
\be
\tilde{d}(\nb,\kb)
\le
d(\nb,\kb).
\en
This implies that in \eqref{ineq} 
the equality always takes place.
Hence $\iota$ and $\phi$ are isomorphisms. 
\qed

\subsection{Characters and $q$-supernomial
coefficients}\label{subsec:2.10}
In this section we identify the characters of the fusion product 
with the $q$-supernomial coefficients for $\sln$ 
known in the literature \cite{HKKOTY}, \cite{SW2}.
Our consideration is restricted to the case of symmetric 
tensor representations. 

For the description of the characters, we find it convenient to 
use the language of partitions. 
For a partition $\lambda=(\lambda_1,\cdots,\lambda_N)$, 
let $\lambda'$ denote the conjugate partition. 
We write $|\lambda|=\sum_{i=1}^N\lambda_i$.
We identify partitions obtained by concatenating $0$'s at the end. 
Reordering $n_p,k_p$ by using \eqref{Wsym}, we can assume without loss of generality that 
\bea
&&n_1\ge \cdots\ge n_N,
\label{order5}\\
&&k_{N_a+1}\ge \cdots\ge k_{N_{a-1}}
\quad (a=1,\cdots,n).
\label{order6}
\ena
Define partitions $\mu^{(1)}\subset \cdots\subset \mu^{(n)}$ by setting 
\be
&&\mu^{(a)}=\sum_{b=1}^a\kappa^{(b)},
\\
&&\kappa^{(b)}=(k_{N_b+1},\ldots,k_{N_{b-1}})'.
\en
Note that $\mu^{(a)}_1=N-N_a$. 
We adopt the following notation
\bea
W^{(n)}[\mu^{(n)},\cdots,\mu^{(1)}]=W(\nb,\kb).
\label{Wmu}
\ena
When $k_1=\cdots=k_{N_{n-1}}=0$, we have $\mu^{(n)}=\mu^{(n-1)}$.
In this case we understand that 
\be
W^{(n)}[\mu^{(n-1)},\mu^{(n-1)},\cdots,\mu^{(1)}]
=
W^{(n-1)}[\mu^{(n-1)},\cdots,\mu^{(1)}].
\en
Assuming $\mu^{(n)}\neq \mu^{(n-1)}$ 
and therefore $n_1=n$, let $k=k_1$.
Theorem \ref{thm:2.2} implies the following exact sequence. 
\bea
&&
\label{recur0}\\
&&
0\longrightarrow 
W^{(n)}[\mu^{(n)},
\mu^{(n-1)}+(1^k),
\cdots,\mu^{(1)}]
\nn
\\
&&\quad 
\overset{\iota}{\longrightarrow}
W^{(n)}[\mu^{(n)},\mu^{(n-1)},\cdots,\mu^{(1)}]
\overset{\psi}{\longrightarrow}
W^{(n)}[\mu^{(n)}-{\bf e}_k,\mu^{(n-1)},\cdots,\mu^{(1)}]
\longrightarrow 0.
\nn
\ena
Here 
${\bf e}_k$ denotes the unit vector $(\delta_{jk})_{1\le j\le k}$ and $(1^k)=\sum_{i=1}^k{\bf e}_i$. 
The subspace corresponds to changing $n_1$ to $n_1-1$, 
while the quotient space corresponds to 
changing $k_1$ to $k_1-1$. 

Let us rewrite \eqref{recur0} in terms of the character (see \eqref{chara})
\bea
&&\chi^{(n)}[\mu^{(n)},\cdots,\mu^{(1)}](q,z_1,\cdots,z_{n-1})
\label{chi}\\
&&\quad
={\rm ch}_{q,z_1,\cdots,z_{n-1}} W^{(n)}[\mu^{(n)},\cdots,\mu^{(1)}]
\nn\\
&&\quad =\sum q^d z_1^{m_1}\cdots z_{n-1}^{m_{n-1}}
\dim W^{(n)}[\mu^{(n)},\cdots,\mu^{(1)}]_{d,m_1,\cdots,m_{n-1}}.
\nn
\ena
\begin{prop}\label{prop:2.5}
We have 
\bea
&&
\chi^{(n)}[\mu^{(n)},\cdots,\mu^{(1)}](q,z_1,\cdots,z_{n-1})
\label{recur}
\\
&&
=
\chi^{(n)}[\mu^{(n)},\mu^{(n-1)}+(1^k),
\cdots,\mu^{(1)}](q,z_1,\cdots,qz_{n-1})
\nn\\
&&
+z_{n-1}\chi^{(n)}[\mu^{(n)}-{\bf e}_k,\cdots,\mu^{(1)}](q,z_1,\cdots,z_{n-1}).
\nn
\ena
\end{prop}
It is useful to modify the character 
\be
&&\tilde\chi^{(n)}
[\mu^{(n)},\cdots,\mu^{(1)}](q,z_1,\cdots,z_{n-1})
\\
&&=
\chi^{(n)}[\mu^{(n)},\cdots,\mu^{(1)}](q,q^{\mu^{(1)}_1}z_1,\cdots, q^{\mu^{(n-1)}_1}z_{n-1}),
\en
so as to absorb 
the shift of variable $z_{n-1}\rightarrow qz_{n-1}$.  
Eq. \eqref{recur} then becomes
\begin{eqnarray}
&&
\label{MODREL}\\
&&\tilde\chi^{(n)}[\mu^{(n)},\cdots,\mu^{(1)}]
\nn\\
&&=\tilde\chi^{(n)}
[\mu^{(n)},\mu^{(n-1)}+(1^k),
\cdots,\mu^{(1)}]
+q^{\mu^{(n-1)}_1}z_{n-1}
\tilde\chi^{(n)}[\mu^{(n)}-{\bf e}_k,\mu^{(n-1)},\cdots,\mu^{(1)}]
\nn
\end{eqnarray}
where we have suppressed the variables 
$(q,z_1,\cdots,z_{n-1})$ which are common to all terms. 
The base of the recursion is 
\begin{equation*}
\tilde\chi^{(1)}[\mu^{(1)}]=1.
\end{equation*}

Using \eqref{recur0} repeatedly, 
we are led to 
\begin{prop}\label{prop:2.6}
There exists a filtration by submodules 
\be
W^{(n)}[\mu^{(n)},\cdots,\mu^{(1)}]=F^0\supset F^1\supset F^2\supset\cdots
\en
such that each composition factor $F^j/F^{j-1}$ has the form 
\be
W^{(n-1)}[\nu^{(n-1)},\mu^{(n-2)},\cdots,\mu^{(1)}],
\en
with appropriate shift of gradings, where $\nu^{(n-1)}$
is a partition satisfying  $\mu^{(n)}\supset \nu^{(n-1)} \supset\mu^{(n-1)}$.
\end{prop}
In terms of the characters, it means a relation of the form
\bea
&&\tilde\chi^{(n)}[\mu^{(n)},\mu^{(n-1)},\cdots,\mu^{(1)}]
\label{rec2}
\\
&&
=\sum_{\mu^{(n)}\supset \nu^{(n-1)} \supset\mu^{(n-1)}}
C_{\nu^{(n-1)}}(q,z_{n-1})
\tilde\chi^{(n-1)}[\nu^{(n-1)},\mu^{(n-2)},\cdots,\mu^{(1)}],
\nn
\ena
where $C_{\nu^{(n-1)}}(q,z_{n-1})$ are some polynomials in $q$ and $z_{n-1}$. 
\medskip

\noindent{\it Example.}\quad 
The following diagram illustrates the process of the 
reduction on the example
$n=2$, $k=2$,
$\mu^{(2)}=(3,2)$, $\mu^{(1)}=\emptyset$ and $\nu^{(1)}=(2,1)$.
{
\setlength{\unitlength}{.5\unitlength}
\begin{eqnarray*}
\begin{CD}
\left(\rb{-15}{\bp(45,30) \row3{30} \row2{15} \ep},
\emptyset\right)
 @>{z_1}>>
\left(\rb{-15}{\bp(45,30) \row3{30} \row1{15} \ep},
\emptyset\right) @. {}
\\
@AAA @AAA @. \\
\left(\rb{-15}{\bp(45,30) \row3{30} \row2{15} \ep},
\rb{-15}{\bp(15,30) \row1{30} \row1{15} \ep}\right)
@>{qz_1}>>
\left(\rb{-15}{\bp(45,30) \row3{30} \row1{15}
    \ep},\rb{-15}{\bp(15,30) \row1{30} \row1{15} \ep}\right) 
@>{qz_1}>>
\left(\rb{-15}{\bp(30,30) \row2{30} \row1{15} \ep},
  \rb{-15}{\bp(15,30) \row1{30} \row1{15} \ep}\right)
\\
@. @AAA @AAA  \\
{} @. \left(\rb{-15}{\bp(45,30) \row3{30} \row1{15}\ep},
\rb{-15} {\bp(30,30)  \row2{30} \row1{15} \ep}
\right)
@>{q^2z_1}>>
\left(\rb{-15} {\bp(30,30)  \row2{30} \row1{15} \ep},\rb{-15} {\bp(30,30)  \row2{30} \row1{15} \ep}\right)
\\
\end{CD}
\end{eqnarray*}
}

The vertical and the horizontal arrows correspond 
to the maps $\iota$ and $\psi$, respectively.
The coefficients appearing in \eqref{MODREL} 
are as indicated.
\qed
\medskip

Let us compute the coefficients $C_{\nu^{(n-1)}}(q,z_{n-1})$.
Assuming $\mu^{(n)}\neq \mu^{(n-1)}$, 
we apply \eqref{rec2} repeatedly for $k_1=k,k-1,\cdots$, 
and collect the terms
which have the array of partitions of the form 
($\nu^{(n-1)}, \nu^{(n-1)},\mu^{(n-2)},\cdots$). 
First consider $k_1=k$. 
If we take the first term from \eqref{rec2} 
$m_1$ times and the second term 
$m_2$ times, then the resulting array of partitions
becomes 
($\mu^{(n)}-m_2\mathbf{e}_k, \mu^{(n-1)}+m_1(1^k),\mu^{(n-2)},\cdots$). 
Comparing the $k$-th row, we have $m_1=\nu^{(n-1)}_k-\mu^{(n-1)}_k$ and
$m_2=\mu^{(n)}_k-\nu^{(n-1)}_k$.  
There are altogether $\displaystyle{\binom{m_1+m_2}{m_1}}$ such terms. 
It is easy to see that the sum of the corresponding coefficients is
\begin{eqnarray*}
(q^{\mu^{(n-1)}_1}z)^{\mu^{(n)}_k-\nu^{(n-1)}_k}
\left[\mu^{(n)}_k-\mu^{(n-1)}_k \atop
\nu^{(n-1)}_k-\mu^{(n-1)}_k \right].
\end{eqnarray*}
One can continue a similar process for $k_1=k-1,\cdots,1$. 
The coefficient $C_{\nu^{(n-1)}}(q,z_{n-1})$ 
is the product of all factors resulting from each step.

Let us s summarize the above calculation.

For a pair of partitions $\kappa\supset\nu$, define 
\bea
F_{\kappa,\nu}(q)=
q^{\sum_{a=1}^{k-1}\nu_{a+1}(\kappa_a-\nu_a)}
\prod_{a=1}^{k}
\left[{\kappa_a-\nu_{a+1}\atop\nu_a-\nu_{a+1}}\right].
\label{tane}
\ena
Here $k=\kappa_1'$ is the length of $\kappa$, 
and the $q$-binomial symbol is 
\be
&&\qbin{m}{n}=\begin{cases}
\displaystyle\frac{[m]!}{[n]![m-n]!} & (0\le n\le m),\\
0 & \mbox{otherwise},\\
\end{cases}
\en
with $[n]!=\prod_{j=1}^n(1-q^j)/(1-q)$.
Returning to the definition \eqref{chi},
we find the following.
\begin{thm}\label{thm:2.4}
We have
\bea
&&
\chi^{(n)}[\mu^{(n)},\mu^{(n-1)},\cdots,\mu^{(1)}](q,z_1,\cdots,z_{n-1})
\label{FERM1}
\\
&&=
\sum_{\mu^{(n)}\supset \nu\supset \mu^{(n-1)}}
z_{n-1}^{|\mu^{(n)}|-|\nu|}
F_{\mu^{(n)}-\mu^{(n-1)},\nu-\mu^{(n-1)}}(q)
\nn
\\
&&\quad \times
\chi^{(n-1)}[\nu,\cdots,\mu^{(1)}](q,z_1,\cdots,z_{n-2}).
\nn
\ena
Hence we obtain a fermionic expression for the character 
\bea
&&
\chi^{(n)}[\mu^{(n)},,\cdots,\mu^{(1)}](q,z_1,\cdots,z_{n-1})
\label{FERM}\\
&&=
\sum\prod_{a=1}^{n-1}\Bigl(
z_a^{|\nu^{(a+1)}|-|\nu^{(a)}|}
F_{\nu^{(a+1)}-\mu^{(a)},\nu^{(a)}-\mu^{(a)}}(q)
\Bigr),
\nn
\ena
where the sum extends over sequences of partitions 
\[
\mu^{(n)}=\nu^{(n)}\supset\nu^{(n-1)}\supset\cdots\supset\nu^{(1)}
\]
satisfying $\nu^{(a)}\supset\mu^{(a)}$ $(1\leq a\leq n-1)$.
\end{thm}

In the case $n_1=\cdots=n_N=n$, we have 
$\mu^{(1)}=\cdots=\mu^{(n-1)}=\emptyset$ and 
$\mu^{(n)}=\mu$ where $\mu=(k_1,\cdots,k_N)'$. 
By comparing \eqref{FERM} with known results
in the literature, we conclude that 
the coefficients of monomials in $z_a$ 
are the $q$-supernomial coefficients. 
More precisely, setting 
\be
W^{(n)}[\mu]=W^{(n)}[\mu,\emptyset,\cdots,\emptyset],
\quad
\chi^{(n)}[\mu]=\chi^{(n)}[\mu,\emptyset,\cdots,\emptyset],
\en
we have 
\bea
\chi^{(n)}[\mu](q,z_1,\cdots,z_{n-1})
=\sum_{\lambda\in\Z^n_{\ge 0}\atop
|\lambda|=|\mu|}
z_1^{\lambda_2}\cdots z_{n-1}^{\lambda_{n}} 
\tilde{S}_{\lambda,\mu}(q).
\label{SS}
\ena
Here $\lambda=(\lambda_1,\cdots,\lambda_n)$ 
and $\tilde{S}_{\lambda,\mu}(q)$ denotes the 
$q$-supernomial coefficient as given in 
formula (2.1) of \cite{S}.  

The exact sequence \eqref{recur0} gives
also a recursive method of constructing 
a monomial basis of $W^{(n)}[\mu^{(n)},\cdots,\mu^{(1)}]$.
\begin{thm}\label{thm:2.3}
There exists a basis $B^{(n)}[\mu^{(n)},\cdots,\mu^{(1)}]$
of $W^{(n)}[\mu^{(n)},\cdots,\mu^{(1)}]$
defined recursively as follows.  
\begin{enumerate}
\item
Each element in the set $B^{(n)}[\mu^{(n)},\cdots,\mu^{(1)}]$
is represented by a monomial in the variables $\{x_b[i]\}$. 
\item $B^{(1)}[\mu^{(1)}]=\{1\}$.
\item If the length of $\mu^{(n)}-\mu^{(n-1)}$ is 
$k>0$, then 
\be
&&B^{(n)}[\mu^{(n)},\cdots,\mu^{(1)}]
\\
&&\quad
=
\iota\left( 
B^{(n)}[\mu^{(n)},\mu^{(n-1)}+(1^k),
\cdots,\mu^{(1)}]\right)
\bigsqcup
x_a[0]
B^{(n)}[\mu^{(n)}-{\bf e}_k,\cdots,\mu^{(1)}],
\en
where the second term in the right hand side 
means elementwise multiplication by $x_a[0]$.  
\end{enumerate}
\end{thm}

\noindent{\it Example.}\quad 
Consider the case of $\slt$ (with $e\overset{\rm def}{=}x_1$). 
Computing the basis for the example
above, with $\mu^{(2)}=(3,2)$ and $\mu^{(1)}=\emptyset$, we have
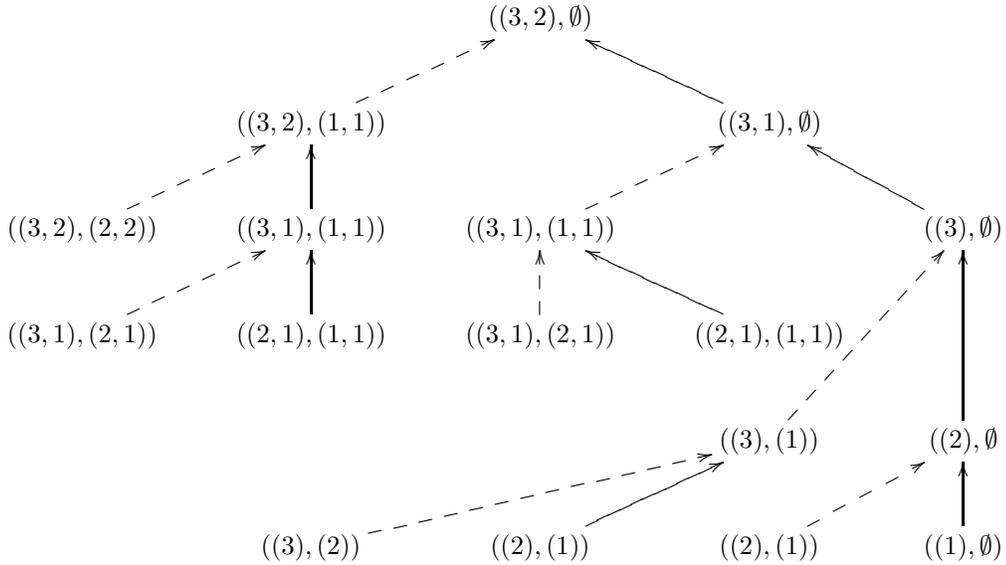
\begin{figure}[here]
{\small
\xymatrix{
{} &{} &{{((3,2),{\emptyset})} }  &&{}&{} &{} &{} &{} &{} &{}
\\
{} &{((3,2),(1,1))}\ar@{-->}[ur] &{} &{((3,1),\emptyset)}\ar[ul] &{} &{} &{} &{} &{} &{} &{} &{} \\
{((3,2),(2,2))}\ar@{-->}[ur]  &{((3,1),(1,1))}\ar[u]&{((3,1),(1,1))}\ar@{-->}[ur] &{}&{((3),\emptyset)}\ar[ul] &{} &{} &{} &{} &{} &{} &{} \\
{((3,1),(2,1))}\ar@{-->}[ur] &{((2,1),(1,1))}\ar[u] 
&{((3,1),(2,1))}\ar@{-->}[u] &
{((2,1),(1,1))}\ar[ul] \\
&{}&{} &{((3),(1))}\ar@{-->}[uur] &{((2),\emptyset}\ar[uu] \\
{}  &{((3),(2))}\ar@{-->}[urr] &{((2),(1))}\ar[ur] &{((2),(1))}\ar@{-->}[ur] &{((1),\emptyset)}\ar[u] &{} &{} &{} &{} &{} \\
}
}
\end{figure}

Here, the dotted arrows correspond to the map $\iota: e[i]\to e[i+1]$,
and the solid arrows to multiplication by $e[0]$. As a result, and
noting that $B^{(2)}[\mu+(1),\mu]=\{1,e[0]\}$, we have
\begin{eqnarray*}
B^{(2)}[(3,2),\emptyset] &=&
\{1, e[0], e[0]^2, e[0]^3, e[0]^4, e[0]^5, e[1],  e[0]e[1],\\ && e[0]^2
e[1], e[0]^3 e[1], e[2], e[1]^2,
 e[0]e[2], \\ && e[0]e[1]^2, e[0]^2 e[2],
e[0]^2 e[1]^2, e[1]e[2], e[1]^3\}.
\end{eqnarray*}

\setcounter{section}{2}
\setcounter{equation}{0}

\section{Character of $\slth$ spaces of coinvariants}\label{sec:3}
In \cite{FJKLM}, we have introduced 
a large class of spaces of coinvariants, 
and established their equivalence to 
certain quotients of filtered tensor products 
of finite-dimensional modules. 
In this section, we restrict to the simplest case 
where the underlying simple Lie algebra is $\g=\slt$, 
and apply the results of the previous section 
to the character of the space of coinvariants. 

\subsection{Spaces of coinvariants}\label{subsec:3.1}
We retain the following notation of \cite{FJKLM}: 
$e,f,h$ are the standard generators of $\slt$. 
$V_l$ is the ($l+1$)-dimensional 
irreducible representation with lowest weight vector $v_l$. 
$L^{(k)}_l(\infty)$ is the level $k$ integrable module 
of the affine Lie algebra $\slth$ `placed at infinity'. 
It has $V_l$ as its subspace of degree $0$
and $x[i]V_l=0$ if $x\in\slt,i<0$.  
By $\mathcal{V}^{(k)}=\oplus_{l=0}^k\Z\cdot[l]$ we denote 
the level $k$ Verlinde algebra for $\slt$ 
with basis $[l]$ ($l=0,\cdots,k$).
We write an element $a\in \mathcal{V}^{(k)}$ as 
$a=\sum_{l=0}^k(a:[l])_k[l]$, where $(a:[l])_k\in\Z$. 
For an $\slt[t]$-module $L$ and a right ideal $Y\subset U(\slt[t])$, we write
$L/Y$ for the space of coinvariants $L/Y L$.  
Let $\sbld=(m_1,\cdots,m_k)$ be an array 
of non-negative integers.
The level-restricted Kostka polynomial for $\slt$ 
is 
\bea
K^{(k)}_{l,\sbld}(q)=
\sum_{\sbf\in\Z_{\ge 0}^k
\atop 2|\sbf|=|\mb|-l}
q^{\sbf\cdot A\sbf+\mathbf{v}\cdot\sbf}
\qbin{A(\sbld-2\sbf)-\mathbf{v}+\sbf}{\sbf},
\label{resKos}
\ena
where $A$ is a $k\times k$ matrix with entries 
$A_{ij}=\min(i,j)$, 
$v_i=(i-k+l)_+$, $|\mathbf{u}|=\sum_{i=1}^k iu_i$,
and $\mathbf{u}\cdot\mathbf{u'}=\sum_{i=1}^ku_i u'_i$.  
In what follows we use also the partition 
\bea
\mu=(k^{m_k}\,\cdots\,1^{m_1})'
\label{part1}
\ena
to label $\sbld$. 
We abuse the notation and write 
$K^{(k)}_{l,\mu}(q)$ for 
$K^{(k)}_{l,\sbld}(q)$.
We recall also the following formula 
in the Verlinde algebra $\mathcal{V}^{(k)}$, 
which will be used later.
\bea
&&[k]^{m_k}\cdots[1]^{m_1}=
\sum_{l=0}^k K^{(k)}_{l,\mu}(1)\cdot [l].
\label{verlinde}
\ena
When we deal with fusion product, we consider also $\slt$ embedded in $\sltr$. 
We write the generators as $e_{23},e_{32},h_{23}$. 
In particular, we use $h_{23}[0]$ for the action on the fusion product, 
and $h[0]$ for that on $L^{(k)}_l(\infty)$. 

Consider the fusion product of irreducible modules 
\be
W^{(2)}[\mu]=
\overbrace{V_k*\cdots * V_k}^{m_k}
* \cdots*
\overbrace{V_1* \cdots * V_1}^{m_1}
\en
taking the tensor product of lowest weight vectors
\be
\overbrace{v_k\otimes\cdots\otimes v_k}^{m_k}
\otimes \cdots\otimes 
\overbrace{v_1\otimes\cdots\otimes v_1}^{m_1}
\en
as cyclic vector. 
The eigenvalue of $h_{23}[0]$ on it is $-|\mu|$. 
Denote by $I^{(k)}_l$ the left ideal of $U(\slt[t])$ 
generated by the elements
\be
h_{23}[0]+l,\quad e_{23}[0],\quad e_{23}[1]^{k-l+1}. 
\en
Let $S(I^{(k)}_l)$ be the right ideal obtained by applying 
the antipode $S$. 
The following result was obtained in 
\cite{FJKLM},Theorem 4.1.  
\begin{prop}\label{prop:3.3}
We have 
\be
{\rm ch}_{q,z_1} W^{(2)}[\mu]/S(I^{(k)}_l)=
z_1^{\frac{|\mu|-l}{2}}K^{(k)}_{l,\mu}(q).
\en
\end{prop}

\subsection{Fusion product of $\sltr$-modules}\label{subsec:3.2}
In this section, 
we consider reducible $\slt$-modules of the form  
\bea
W_r=V_0\oplus\cdots\oplus V_r. 
\label{Wr}
\ena
It is generated over $\slt$ by 
the sum of lowest weight vectors $w_r=v_0+\cdots+v_r$. 
We make use of the fact that \eqref{Wr} 
can also be viewed as the $r$-th symmetric tensor 
$S^r(V_{(1,0)})$ of the defining representation of $\sltr$, 
restricted to the subalgebra $\slt\subset \sltr$.

Let $\mathbf{M}=(M_1,\cdots,M_k)\in\Z_{\ge 0}^k$, and 
\bea
\lambda=(k^{M_k},\cdots,1^{M_1})'.
\label{part2}
\ena
Fix a set of distinct complex numbers 
$\zz=(z_1,\cdots,z_N)$, $N=M_1+\cdots+M_k$. 
Consider the filtered tensor product of $\slt[t]$-modules 
\bea
\mathcal{F}_{\zz}
\Bigl(\overbrace{W_k,\cdots , W_k}^{M_k}
, \cdots,
\overbrace{W_1, \cdots, W_1}^{M_1}\Bigr),
\label{W2}
\ena
taking 
\bea
\overbrace{w_k\otimes\cdots\otimes w_k}^{M_k}
\otimes \cdots\otimes
\overbrace{w_1\otimes \cdots\otimes w_1}^{M_1}
\label{wM}
\ena
as cyclic vector.
\begin{lem}\label{lem:3.1}
The filtered tensor product \eqref{W2}
as $\slt[t]$-module and with cyclic vector \eqref{wM} 
coincides with the filtered tensor product 
as $\sltr[t]$-module and with cyclic vector 
\bea
\overbrace{v_k\otimes\cdots\otimes v_k}^{M_k}
\otimes \cdots\otimes
\overbrace{v_1\otimes \cdots\otimes v_1}^{M_1}.
\label{vM}
\ena
\end{lem}
The proof is given in Appendix, Proposition \ref{prop:app2}.
Note that the fusion product in the latter sense 
is the same as the one over the abelian Lie algebra
$\ga[t]$, $\ga=\C e_{13}\oplus \C e_{23}$.
That is, 
\bea
W^{(3)}[\lambda]
=
\overbrace{W_k*\cdots * W_k}^{M_k}
* \cdots*
\overbrace{W_1* \cdots * W_1}^{M_1}
\label{W3}
\ena
in the notation of the previous section. 

Let $Y_r$ denote the right 
ideal of $U(\slt[t])$ generated by
\be
e_{23}[0]^{r+1},~~e_{32}[0],~~\quad x[i]\quad (x\in \slt
\subset \mathfrak{sl}_3, i> 0).
\en
As a left $\slt[t]$-module we have 
\be
U(\slt[t])/S(Y_r)\simeq W_r,  
\en
where 
$x[i]$ ($i>0$) acts as $0$ on the right hand side. 
Introduce the fusion right ideal in the sense of \cite{FJKLM}
\bea
Y_{\mub}=
\overbrace{Y_k\fusn \cdots\fusn Y_k}^{M_k}
\fusn
\cdots\fusn\overbrace{Y_1\fusn \cdots\fusn Y_1}^{M_1}(\zz). 
\label{right}
\ena
Our aim is to determine the character of the space of coinvariants
\bea
\mathop{\rm gr}L^{(k)}_{l}(\infty)/Y_{\mub}.
\label{coinv1}
\ena

Concerning \eqref{coinv1}, the following are known. 
\begin{prop}\label{prop:3.1}$($\cite{FJKLM}, Theorem 2.9$)$
The dimension of \eqref{coinv1} is given by the Verlinde rule 
\be
\dim L^{(k)}_l(\infty)/Y_{\mub}=
\Bigl(([0]+\cdots+[k])^{M_k}\cdots([0]+[1])^{M_1}:[l])_k.
\en
\end{prop}
{}From \eqref{verlinde}, the right hand side is written explicitly as 
\bea
&&\sum_{\mu\subset \lambda
\atop |\mu|\equiv l\bmod 2}
\prod_{i=1}^k
\binom{\lambda_i-\mu_{i+1}}{\mu_i-\mu_{i+1}}
\cdot K^{(k)}_{l,\mu}(1)
=\sum_{\mu\subset \lambda
\atop |\mu|\equiv l\bmod 2}
F_{\lambda,\mu}(1)K^{(k)}_{l,\mu}(1).
\label{au1}
\ena

\begin{prop}\label{prop:3.2}$($\cite{FJKLM}, Theorem 3.6$)$
There is an isomorphism of filtered vector space
\bea
&&L^{(k)}_l(\infty)/Y_{\mub}
\simeq
\mathcal{F}_{\zz}(
\overbrace{W_k,\cdots, W_k}^{M_k},\cdots,
\overbrace{W_1,\cdots, W_1}^{M_1})
/S(I^{(k)}_l).
\label{paper1}
\ena
\end{prop}

On the space $L^{(k)}_l(\infty)$ we have a grading by $d=td/dt$ and $h[0]$. 
Define the character of the associated graded space  
$L=\gr L^{(k)}_l(\infty)/Y_\mub$ by 
\be
{\rm ch}_{q,z}L
={\rm tr}_{L}\bigl(q^dz^{h[0]}\bigr).
\en
We have also a grading on the fusion product $W^{(3)}[\lambda]$. 
Its character is defined in \eqref{chara}, by
assigning  $q^iz_1$ (resp. $q^iz_2$) to $e_{23}[i]$ (resp.$e_{13}[i]$)
and $1$ to the cyclic vector. 
Using the action of $h_{ab}[0]\in \sltr$, we can write the character as
\be
\chi^{(3)}[\lambda](q,z_1,z_2)=
(z_1z_2)^{|\lambda|/3}
{\rm tr}_{W^{(3)}[\lambda]}\bigl(
q^d z_1^{\overline{\varepsilon}_2}
z_2^{\overline{\varepsilon}_1}\bigr),
\en
where 
\be
&&\overline{\varepsilon}_1=\frac{1}{3}(h_{12}[0]+h_{13}[0]),
\quad
\overline{\varepsilon}_2=\frac{1}{3}(h_{21}[0]+h_{23}[0]).
\en

\begin{lem}\label{lem:3.2}
In the isomorphism \eqref{paper1}, 
the grading by $h[0]$ in the left hand side translates 
in the right hand side to the grading by the operator 
$(h_{21}[0]+h_{31}[0]+2|\lambda|)/3$. 
\end{lem}
\begin{proof}
Consider the cyclic vector $w_r=v_0+\cdots+v_r$ of $W_r$. 
As it is explained in (3.20), \cite{FJKLM}, 
the grading by $h[0]$ counts $xv_j$ as $j$ for any $x\in \slt$.  
Denote the natural basis of $\C^3$ by $u_i$ ($i=1,2,3$), 
on which $e_{jj}[0]$ acts as $e_{jj}[0]u_i=\delta_{ij}u_i$. 
In the identification $W_r\simeq S^r(\C^3)$, 
$v_j$ is identified with the symmetrization of 
$u_1^{\otimes r-j}\otimes u_3^{\otimes j}$. 
Therefore $h[0]$ acts on $S^r(\C^3)$ as
\be
e_{22}[0]+e_{33}[0]=
\frac{1}{3}\bigl(h_{21}[0]+h_{31}[0]+2r\bigr).
\en
The lemma follows from this.
\end{proof}

{}From the lemma we find the following relation 
\be
{\rm ch}_{q,z}\gr L^{(k)}_l(\infty)/Y_\mub
=z^{|\lambda|}{\rm ch}_{q,z_1,z_2}
W^{(3)}[\lambda]/S(I^{(k)}_l)\bigl|_{z_1=1,z_2=z^{-1}},
\en
where we have used 
\be
z_1^{\overline{\varepsilon}_2}
z_2^{\overline{\varepsilon}_1}
=z_1^{h_{23}[0]/2}(z_1^{1/2}z_2^{-1})^{-h[0]+2|\lambda|/3}.
\en

\subsection{Character formulas}\label{subsec:3.3}
With the preparations above, 
let us proceed to the character of a quotient of the fusion product. 

\begin{prop}\label{prop:3.4}
We have 
\be
&&{\rm ch}_{q,z_1,z_2} W^{(3)}[\lambda]/S(I^{(k)}_l)
=
\sum_{\mu\subset \lambda
\atop |\mu|\equiv l\bmod 2}
z_1^{\frac{|\mu|-l}{2}}z_2^{|\lambda|-|\mu|}
F_{\lambda,\mu}(q)K^{(k)}_{l,\mu}(q),
\en
where $F_{\lambda,\mu}(q)$ is defined in \eqref{tane}.
\end{prop}
\begin{proof}
By Proposition \ref{prop:3.1} and \ref{prop:3.2}, 
\eqref{au1} gives a lower bound to the dimension 
$\dim W^{(3)}[\lambda]/S(I^{(k)}_l)$.

On the other hand, from Proposition \ref{prop:2.6} 
we have a filtration of $W^{(3)}[\lambda]$ 
by submodules with composition factors of the form 
$W^{(2)}[\mu]$. 
The operators $h_{23}[0]+l$, $e_{23}[0]$ and
$e_{23}[1]^{k-l+1}$ leave the filtration invariant. 
Replacing the character of the quotient 
by that of the quotient of the associated graded space, 
we obtain an upper estimate
\bea
&&{\rm ch}_{q,z_1,z_2} W^{(3)}[\lambda]/S(I^{(k)}_l)
\nn\\
&&\quad 
\le
\sum_{\mu\subset\lambda
\atop |\mu|\equiv l\bmod 2}
z_2^{|\lambda|-|\mu|}
F_{\lambda,\mu}(q)
{\rm ch}_{q,z_1} W^{(2)}[\mu]/S(I^{(k)}_l).
\label{au2}
\ena
Proposition \ref{prop:3.3} states that up to a power of $z_1$ 
the last factor is given by $K^{(k)}_{l,\mu}(q)$. 
Setting $q=z_1=z_2=1$, the right hand side of 
\eqref{au2} becomes the lower bound \eqref{au1}.
This in turn implies that the equality holds in each step,  
and the proof is over.
\end{proof}

Writing $F_{\lambda\mu}(q)$ explicitly, we obtain 
the following as an immediate consequence of Proposition \ref{prop:3.3}
and \ref{prop:3.2}. 
\begin{thm}\label{thm:3.1}
We have a fermionic formula  for the character 
of the space of coinvariants 
\bea
&&{\rm ch}_{q,z} \gr L^{(k)}_{l}(\infty)/Y_{\mub}
\label{fer2}\\
&&
\quad=
\sum_{\mu\subset \lambda
\atop |\mu|\equiv l\bmod 2}
z^{|\mu|}
q^{\sum_{i=1}^{k-1}\mu_{i+1}(\lambda_i-\mu_i)}
\prod_{i=1}^k
\qbin{\lambda_i-\mu_{i+1}}{\mu_i-\mu_{i+1}}
\cdot 
K^{(k)}_{l,\mu}(q).
\nn
\ena
\end{thm}

Since $W^{(3)}[\lambda]$ has an $\slt$-structure, 
its character is a linear combination of 
$\slt$-characters 
$\chi_j(z)=z^{-j/2}+\cdots+z^{j/2-1}+z^{j/2}$. 
Let us define $X_{j\lambda}(q,z)$ as the coefficient of the 
expansion
\be
\chi^{(3)}[\lambda](q,z_1,z_2)
=z_2^{|\lambda|}\sum_{j\ge 0}\chi_j(z_1)
X_{j\lambda}(q,z^{1/2}_1z_2^{-1}).
\en
Then we have the following alternating sum formula.
When specialized to the case $\mathbf{M}=(0,\cdots,0,M)$, 
it settles the `main conjecture' of  
\cite{FL} (Conjecture 3.5) for $\mathfrak{g}=\slt$. 

\begin{thm}\label{thm:3.2}
\bea
{\rm ch}_{q,z} \gr L^{(k)}_{l}(\infty)/Y_{\mub}
&=&\sum_{i\ge 0}q^{(k+2)i^2+(l+1)i}
X_{2(k+2)i+l,\lambda}(q,z) 
\label{bos2}\\
&&-\sum_{i> 0}q^{(k+2)i^2-(l+1)i}
X_{2(k+2)i-l-2,\lambda}(q,z).  
\nn
\ena
\end{thm}
\begin{proof}
We have from \eqref{FERM1} and \eqref{SS} 
\bea
\chi^{(3)}[\lambda](q,z_1,z_2)
&=&\sum_{\mu\subset\lambda}z_2^{|\lambda|-|\mu|}F_{\lambda\mu}(q)
\chi^{(2)}[\mu](q,z_1),
\label{xx1}
\\
\chi^{(2)}[\mu](q,z_1)
&=&\sum_{0\le j\le |\mu|}
z_1^{(|\mu|-j)/2}\widetilde{S}_{\sigma(j),\mu}(q),
\label{xx2}
\ena
where we have set 
$\sigma(j)=\bigl(\frac{|\mu|+j}{2},\frac{|\mu|-j}{2}\bigr)$. 
Let $K_{j,\mu}(q)$ 
denote the non-restricted Kostka polynomial for $\slt$, 
which is related to the $q$-supernomial coefficients 
$\widetilde{S}_{\sigma(j),\mu}(q)$ via 
$\chi^{(2)}[\mu](q,z_1)$ as 
\be
\chi^{(2)}[\mu](q,z_1)
=\sum_{0\le j\le |\mu|}z_1^{|\mu|/2}\chi_j(z_1)K_{j,\mu}(q).
\en
{}From this and \eqref{xx1},\eqref{xx2}, we find an expression
\be
X_{j\lambda}(q,z)=\sum_{|\mu|=j}
z^{|\mu|}F_{\lambda\mu}(q)K_{j\mu}(q).
\en
The statement of Theorem now follows from Theorem \ref{thm:3.1}
together with the known formula (see \cite{SS}, eq.(6.8))
\be
K^{(k)}_{l,\mu}(q) 
&=&\sum_{i\ge 0}q^{(k+2)i^2+(l+1)i}
K_{2(k+2)i+l,\mu}(q) 
\nn\\
&&-\sum_{i> 0}q^{(k+2)i^2-(l+1)i}
K_{2(k+2)i-l-2,\mu}(q). 
\label{alter}
\en
\end{proof}

\appendix
\setcounter{equation}{0}
\section{A Lemma on changing cyclic vectors}
In this appendix we prove one of the steps
used in the main text in a more general setting,
i.e., to show the equality of certain filtered tensor products
of $\slnm[t]$ and $\sln[t]$ representations.
In this paper, we use only the case
of $(\slt,\sltr)$ pair and for symmetric tensor representations.

Since the essential part of the argument is contained in the simplest case
$n=2$, we first formulate it and give a proof to that particular case.
Consider the subalgebra $\C h\otimes\C[t] \subset \slt[t]$.

\begin{prop}\label{APP1}
Let $v_a$ $(1\leq a\leq N)$ be the lowest weight vectors
of finite dimensional irreducible representations of $\slt$.
We denote by $U_a$ the representation generated by $v_a$.
Let $\mathcal{Z}=(z_1,\ldots,z_N)$ be (not necessarily distinct)
complex numbers.
Consider the tensor product of the evaluation representations $U_a$
$(1\leq a\leq N)$ of $\slt[t]$, for which the evaluation
parameter is set to $z_a$. We denote this tensor product by $W(\mathcal{Z})$.

Let $w_0=v_1\otimes\cdots\otimes v_N$, and consider the subspace
$G^i$ of $W(\mathcal{Z})$ spanned by the vectors
$x_1[i_1]\cdots x_l[i_l]w_0$ where $l\geq0$,
$x_1,\ldots,x_l\in \slt$ and $\sum_ai_a\leq i$.

Let
\[
w=(\sum_i\frac{e[0]^i}{i!}v_1)\otimes\cdots\otimes(\sum_i\frac{e[0]^i}{i!}v_N),
\]
and let $F^i$ be the subspace of $W(\mathcal{Z})$ spanned by the vectors
$h[i_1]\cdots h[i_l]w$ where $l\geq0$ and $\sum_ai_a\leq i$.

Under this setting, the equality $F^i=G^i$ holds for all $i\geq0$.
\end{prop}
\begin{proof}
Note that
\[
w=\sum_i\frac{e[0]^i}{i!}w_0.
\]
Since the $h[0]$-weights of the vectors $e[0]^iw_0$ $(i\geq0)$ are distinct,
they all belong to $F^0$. Since $\sum_i{\bf C}e[0]^iw_0$ is
$h\otimes t^0$-invariant, we have $F^0=\sum_i{\bf C}e[0]^iw_0$.
Since $f[0]w_0,h[0]w_0\in{\bf C}w_0$, $F^0$ is $\slt\otimes t^0$ invariant,
and therefore, $F^0=G^0$.

Now, we will show that $(\slt\otimes t^j)F^i\subset F^{i+j}$
for all $i,j\geq0$. From this follows $F^i=G^i$ for all $i$.
Note that
\[
F^i=\sum_{m\geq0}\sum_{i_1+\cdots+i_K\leq i}
{\bf C}h[i_1]\cdots h[i_K]e[0]^mw_0.
\]

First, we show that
\[
e[j]h[i_1]\cdots h[i_K]e[0]^mw_0\in F^{i_1+\cdots i_K+j}
\]
by induction on $K$. For $K=0$, we have
\begin{eqnarray*}
e[j]e[0]^mw_0&=&e[0]^me[j]w_0\\
&=&\frac12e[0]^m[h[j],e[0]]w_0.
\end{eqnarray*}
We use $a\equiv_ib$ meaning $a-b\in F^i$.
Since $h[j]w_0\in{\bf C}w_0$ and $h[j]e[0]^{m+1}w_0\in F^j$, we have
\begin{eqnarray*}
e[j]e[0]^mw_0&\equiv_j&\frac12e[0]^mh[j]e_0w_0\\
&\equiv_j&\frac12[e[0]^m,h[j]]e[0]w_0\\
&=&-me[j]e[0]^mw_0.
\end{eqnarray*}
{}From this follows that $e[j]e[0]^mw_0\in F^j$.

For $K\geq1$, we have
\[
e[j]h[i_1]\cdots h[i_K]e[0]^mw_0
=[e[j],h[i_1]]\cdots h[i_K]w_0+
h[i_1]e[j]\cdots h[i_K]w_0.
\]
Using induction hypothesis, we see that this element belongs
to $F^{i_1+\cdots i_K+j}$.
\end{proof}

Now, we proceed to the case of $(\slnm,\sln)$ pair.
\begin{prop}\label{prop:app2}
Let $v_a$ $(1\leq a\leq N)$ be the lowest weight vectors
of finite dimensional irreducible representations of $\sln$ $(n\geq3)$.
We denote by $U_a$ the representation generated by $v_a$. Let
$\mathcal{Z}=(z_1,\ldots,z_N)$ be (not necessarily distinct) complex numbers.
Consider the tensor product of the evaluation representations $U_a$
$(1\leq a\leq N)$ of $\sln[t]$, for which the evaluation
parameter is set to $z_a$. We denote this tensor product by $W(\mathcal{Z})$.

Let $w_0=v_1\otimes\cdots\otimes v_N$, and consider the subspace
$G^i$ of $W(\mathcal{Z})$ spanned by the vectors
$x_1[i_1]\cdots x_l[i_l]w_0$ where $l\geq0$,
$x_1,\ldots,x_l\in \sln$ and $\sum_ai_a\leq i$.

Let
\[
w=(\sum_i\frac{e_{1n}^i}{i!}v_1)\otimes\cdots\otimes
(\sum_i\frac{e_{1n}^i}{i!}v_N),
\]
and let $F^i$ be the subspace of $W(\mathcal{Z})$ spanned by the vectors
$x_1[i_1]\cdots x_l[i_l]w$ where $l\geq0$, $x_1,\ldots,x_l\in \slnm$
and $\sum_ai_a\leq i$.

Under this setting, the equality $F^i=G^i$ holds for all $i\geq0$.
\end{prop}
\begin{proof}
{\it First Step}\,: We show that $F^0$ is spanned by the vectors
\begin{equation}
e_{a_1b_1}[0]\cdots e_{a_Jb_J}[0]\prod_{2\leq a\leq n}e_{1a}[0]^{m_a}w_0,
\label{VEC}
\end{equation}
where $J\geq0$, $2\leq a_i<b_i\leq n$ for $1\leq i\leq J$,
and $m_a\geq0$ for $2\leq a\leq n$.

Proceeding as in the proof of Proposition \ref{APP1}, we see that
$F^0$ contains the vectors $e_{1n}[0]^mw_0$ $(m\geq0)$.
Then, by using $[e_{na}[0],e_{1b}[0]]=-\delta_{n,b}e_{1a}[0]$
and $[e_{1a}[0],e_{1b}[0]]=0$ for $2\leq a\leq n-1$ and $2\leq b\leq n$,
we see that $F^0$ contains the vectors
$\prod_{2\leq a\leq n}e_{1a}[0]^{m_a}w_0$ for $m_2,\ldots,m_n\geq0$.

Now, by using the PBW theorem, it is straightforward to see that
$F^0$ is spanned by the vectors (\ref{VEC}).

{\it Second Step}\,: We show that $F^0=G^0$.
It is enough to show that $F^0$ is invariant by the actions of
$e_{1a}[0]$ and $e_{a1}[0]$ for $2\leq a\leq n$. The invariance for
$e_{1a}[0]$ follows from $[e_{1a}[0],e_{bc}[0]]=\delta_{ab}e_{1c}[0]$
for $2\leq a\leq n$ and $2\leq b<c\leq n$. The invariance for $e_{a1}[0]$
follows from $e_{a1}[0]w_0=0$, $[e_{a1}[0],e_{bc}[0]]=-\delta_{ac}e_{b1}[0]$
for $2\leq a\leq n$ and $2\leq b<c\leq n$,
and
\begin{eqnarray*}
[e_{a1}[0],e_{1b}[0]]=
\begin{cases}
-h_{1a}[0]&\hbox{if }a=b;\\
e_{ab}[0]&\hbox{if }a\not=b,
\end{cases}
\end{eqnarray*}
for $2\leq a,b\leq n$.

{\it Third Step}\,: We show that $F^j$ contains the vectors
\begin{equation}
\label{VECj}
e_{a_1b_1}[j_1]\cdots e_{a_Jb_J}[j_J]\prod_{i=1}^Kh_{c_i}[k_i]
\prod_{i=1}^Le_{1d_i}[l_i]e_{1n}[0]^mw_0,
\end{equation}
where $J,K,L,j_i,l_i,m\geq0$, $k_i>0$,
$2\leq a_i<b_i\leq n$, $2\leq c_i,d_i\leq n-1$, and
\[
\sum_{i=1}^Jj_i+\sum_{i=1}^Kk_i
+\sum_{i=1}^Ll_i\leq j.
\]

We have shown that $F^0$ contains the vectors
$\prod_{a=2}^ne_{1a}[0]^{m_a}w_0$.
Applying $e_{32}[l^{(2)}_i]$ $(1\leq i\leq m'_2\leq m_3;l^{(2)}_i>0)$,
we get
\[
\prod_{a=4}^ne_{1a}[0]^{m_a}e_{13}[0]^{m_3-m'_2}
\prod_{i=1}^{m'_2}e_{12}[l^{(2)}_i]
e_{12}[0]^{m_2}w_0\in F^{\sum_{i=1}^{m'_2}l^{(2)}_i}
\]
Repeating this procedure successively for $e_{43},\ldots,e_{n,n-1}$,
we obtain the vectors
$\prod_{i=1}^Le_{1d_i}[l_i]e_{1n}[0]^mw_0\in F^{\sum_{i=1}^Ll_i}$.
Therefore, $F^j$ contains the vectors (\ref{VECj}).
We denote by $H^j$ the linear span of these vectors.

{\it Last Step}\,: We show that
\[
(\sln\otimes t^i)H^j\subset H^{i+j}.
\]
Since $H^0=F^0=G^0$, it then follows that $H^j=F^j=G^j$.
The only non-trivial calculation is to show that $e_{1n}[i]H^j\subset H^{i+j}$.
This can be shown similarly as in the proof of Proposition \ref{APP1}.
\end{proof}
\bigskip

\def \b {{\mathfrak b}}

Let us explain these statements in more invariant terms. To define a
lowest vector of an irreducible representation it is enough to choose
a Borel subalgebra $\b \subset \sln$. The set of all possible Borel
subalgebras inherits the natural topology of a homogeneous space,
namely they form the flag manifold. 

\def \ou {{\overline{U}}}

Now take the standard inclusion $\slnm \subset \sln$ and choose a
Borel subalgebra $\b \subset \sln$. Then for any irreducible
representation $U$ of $\sln$ by $v_U^\b$ denote the corresponding
lowest weight vector and by $\ou$ denote the restriction of $U$ to
$\slnm$. 

\begin{thm}
Consider a set of irreducible
representations $U_1, \dots, U_N$ of $\sln$ and a set of distinct points
$\zz = (z_1, \dots, z_n)$. Choose a generic Borel subalgebra $\b \subset
\sln$. Then vectors  $v_{U_a}^\b$ are cyclic in
$\ou_a$ as well as in $U_a$ and
we have an isomorphism of filtered spaces
\begin{equation}\label{tres}
{\mathcal F}_\zz (U_1, \dots, U_N) \cong {\mathcal F}_\zz (\ou_1,
\dots, \ou_N).
\end{equation} 
\end{thm}

\begin{proof}
Note, that $\eqref{tres}$ is satisfied for an algebraically open set of
Borel subalgebras. So it is enough to show it for a certain
subalgebra.  

Let $\b_0$ be the standard Borel subalgebra of upper-triangular
matrices. Then Proposition \ref{prop:app2} implies \eqref{tres} for
$$\b = \exp(e_{1n}) \b_0 \exp(-e_{1n}).$$
\end{proof}

\noindent
{\it Acknowledgments.}\quad
This work is partially supported by
the Grant-in-Aid for Scientific Research (B2)
no.12440039, no.14340040  
and (A1) no.13304010, Japan Society for 
the Promotion of Science.
BF is partially supported by grants RFBR 02-01-01015 and
INTAS-00-00055.
The work of SL is partially
supported by the grant RFBR-01-01-00546.
The last stage of this work was carried out
while the authors were visiting Mathematical Sciences
Research Institute, Berkeley, March 2002.

\end{document}